\documentclass[11pt,a4paper]{amsart}
\usepackage{amsmath, amssymb}
\usepackage{array}
\usepackage{longtable}
\usepackage{epstopdf}
\usepackage{graphicx}
\usepackage{xcolor}
\newcommand{\colored}[1]{{\color{black}#1}}

\newtheorem{Thm}{Theorem}[section]
\newtheorem{Prop}[Thm]{Proposition}
\newtheorem{Lem}[Thm]{Lemma}

\newtheorem{Thm*}{Theorem}[section]
\newtheorem{Def}[Thm]{Definition}
\newtheorem{Rem}[Thm]{Remark}

\newcommand{\X}{X \cup X^{-1}} 
\newcommand{\comp}{\mathcal{C}} 
\newcommand{\compb}{\mathcal{D}} 
\newcommand{\graph}{\Gamma} 
\def\M#1{M(#1)} 
\def\init#1{\alpha(#1)} 
\def\ter#1{\omega(#1)} 
\def\pcat#1{\operatorname{IC}(#1)} 
\def\cat{C} 

\newcommand{\Stab}{\operatorname{Stab}}
\newcommand{\FIM}{\operatorname{FIM}}
\newcommand{\SIM}{\operatorname{SIM}}
\newcommand{\bd}{\partial}
\newcommand{\bl}{bl}
\newcommand{\ol}{\overline}

\renewcommand{\int}{\operatorname{int}}
\renewcommand{\star}{\operatorname{star}}

\title{Inverse monoids and immersions of  $\Delta$-complexes}

\author{John Meakin}
\email{jmeakin@math.unl.edu}
\address{Department of Mathematics,
	University of Nebraska-Lincoln\\
	Lincoln, NE 68588, USA}

\author{N\'ora Szak\'acs }
\email{szakacsn@math.u-szeged.hu}
\address{Bolyai Institute,
	University of Szeged\\
	Aradi v\'ertan\'uk tere 1.
	H-6720 Szeged, Hungary}
\address{Center of Mathematics, University of Porto\\ Rua do Campo Alegre 687; 4169-007 Porto, Portugal}
\thanks{The second author was partially supported by the Hungarian National Foundation
	for Scientific Research grant no. K115518; furthermore by CMUP (UID/MAT/00144/2013), which is funded by FCT (Portugal) with national (MEC) and European structural funds (FEDER), under the partnership agreement PT2020.}


\begin{document}

\begin{abstract}

\colored{An immersion $f : {\mathcal D} \rightarrow \mathcal C$ between $\Delta$-complexes is a $\Delta$-map that
induces injections from star sets of $\mathcal D$ to star sets of $\mathcal C$}.
We study immersions
between finite-dimensional connected $\Delta$-complexes by replacing the
fundamental group of the base space by an appropriate inverse
monoid. We show how conjugacy classes of the closed inverse
submonoids of this inverse monoid may be used to classify connected
immersions into the complex. This extends earlier results of
Margolis and Meakin for immersions between graphs and of Meakin and
Szak\'acs on immersions into $2$-dimensional $CW$-complexes.

\end{abstract}

\maketitle

\section{Introduction}

The notion of immersion arises from differential geometry: it is a
differentiable function between differentiable manifolds whose
derivative is everywhere injective. An immersion is essentially a
local smooth embedding: a typical example is the immersion of the
Klein bottle into $3$-space --- it is not an embedding, but it is a
local embedding, which suffices for the purpose of visualization.

In the absence of a differentiable structure, one can define a
topological notion of immersion called a {\em topological immersion}, that
is, a continuous map which is a local homeomorphism onto its image.
Every immersion between differentiable structures is a topological immersion. In our paper, we consider
a combinatorial notion of
 immersions between connected, finite-dimensional
$\Delta$-complexes, which are, as usual in this setting, also
assumed to take $k$-cells to $k$-cells and commute with the
characteristic maps. This definition is the higher-dimensional generalization of graph immersions in the sense of Stallings \cite{stall}, which have played a
significant role in the
study of subgroups of free groups.
Every topological immersion is an immersion in this sense
and for $\Delta$-maps between locally compact $\Delta$-complexes, the converse is true, which we prove in the Appendix.

Covering maps are also immersions, and unlike immersions,  for ``nice
enough'' topological spaces, they are very well understood by means of
the fundamental group of the base space. This classification
makes use of the fact that every path in the base space lifts uniquely
to any given point in the preimage of its initial point; furthermore,
homotopic paths lift to homotopic paths.
It follows that in order to characterize covering maps over a space,
it is sufficient to know which closed paths lift to closed paths under
the covering map, and this is encoded by the fundamental groups of the
two spaces.


In the case of immersions, there is no unique path lifting in the
sense described above; there is however a unique partial path lifting,
that is, paths  lift partially, and this partial lifting is unique at
every point. Furthermore, even when they do lift, homotopic
paths do not necessarily lift to homotopic paths, so homotopy is not
the correct equivalence relation to encode immersions. In order to
obtain a similar characterization to the classification of covers, we
need an algebraic structure rich enough to differentiate between
non-homotopic paths, and to encode not just when closed paths lift to
closed paths,  but also when paths lift at all. This seems difficult
in general, but manageable at least for certain classes of cell
complexes that have a strong enough combinatorial structure, in
particular for $\Delta$-complexes and for immersions between these
complexes that respect the cell structure.

In the cell complex setting, there is a unique partial lifting of
cells of arbitrary dimension and there are complex interrelationships
between liftings of cells of different dimensions. We need to encode
this information algebraically. For cell complexes of dimensions $1$
and $2$ this has been accomplished by introducing an \textit{inverse
monoid} that serves the role of the fundamental group (see the papers
\cite{MM1} for graphs and \cite{MSz} for $2$-dimensional complexes).
However, in the case of higher dimensional cell complexes this
requires significant additional care.

In this paper, we introduce the notion of a \textit{generalized path},
which is essentially a sequence of connecting cells of any dimension.
We then define an equivalence relation on generalized paths that
serves the role of homotopy
in the theory of immersions.
We construct an
inverse monoid of equivalence classes of generalized paths  which we
call a \textit{loop monoid}, and prove that loop monoids encode
immersions between $\Delta$-complexes in the way the fundamental group
encodes coverings.

Loop monoids are defined in Section 3: the loop monoid of a
complex $\comp$ at a point $v$ is denoted by $L(\comp, v)$. Section
6 describes the one-to-one correspondence between conjugacy classes of closed
inverse submonoids of loop monoids and immersions into the complex. 
Theorem
\ref{main} shows that immersions give rise to closed inverse submonoids, and conversely, closed inverse submonoids give rise to immersions, while Theorem \ref{conjugate} shows that two closed inverse submonoids of the loop monoid give rise to the same immersion if and only if they are conjugate.

\section{Preliminaries}

 In this section we introduce the
basic notions of inverse monoids and immersions between $\Delta$-complexes that will be used in the remainder of the paper.

\subsection{Inverse monoids}

An {\it inverse monoid} is a monoid $M$ with the property that for
each $a \in M$ there is a unique element $a^{-1}$ (the inverse of
$a$) in $M$ such that

\begin{center}

$a = aa^{-1}a$ and $a^{-1} = a^{-1}aa^{-1}$.

\end{center}

Inverse monoids arise naturally in the study of partial symmetry in
mathematics in much the same way as groups arise in the study of
symmetry. In fact the {\it Wagner-Preston Theorem} states that every
inverse monoid embeds in an appropriate {\it symmetric inverse
monoid} $\SIM(Q)$, i.e. the monoid of all bijections between subsets
of the set $Q$ under the usual composition of partial maps. For this
and many additional  properties of inverse monoids and their
connections with other fields of mathematics, we refer the reader to
the book of Lawson \cite{law}. Some of the most basic properties of
inverse monoids that we will need are listed in the following
proposition.

\begin{Prop}
\label{invprops}
Let $M$ be an inverse monoid with set $E(M) = \{e \in M : e = e^2\}$
of idempotents. Then $E(M)$ is nonempty, and the following hold:

\begin{itemize}

\item{The idempotents of $M$ commute, i.e. $ef = fe$ for all $e,f
\in E(M)$. Thus the set $E(M)$ of idempotents forms a lower
semilattice with respect to $e \wedge f = ef$. In particular, $g \leq ef$ if and only if $g \leq e$ and $g \leq f$ for any $e,f,g \in E(M)$.}

\item{The relation defined on $M$ by $a \leq b$ iff $a = eb$ for
some $e \in E(M)$ is a partial order on $M$, called the {\it natural
partial order} on $M$. The natural partial order is compatibe with
the multiplication and inversion operations in $M$.}


\item{Let $\sigma_{M}$ be a relation on $M$ such that $a \, {\sigma}_{M} \, b$ iff there exists
$c \in M$ with $c \leq a$ and $c \leq b$.  This is a congruence on
$M$, called the {\it minimum group congruence on $M$}. The quotient
$M/{\sigma}_{M}$ is a group, the {\it greatest group homomorphic
image of $M$}. Equivalently, $\sigma_M$ is generated as a congruence by pairs of the form $(aa^{-1}, 1)$, where $a \in M$.}

\end{itemize}

\end{Prop}

Inverse monoids also arise naturally as transition monoids of {\it
inverse automata}, which are automata whose underlying graphs are
edge labeled over an alphabet $\X$ in the sense described below.

Let $X$ be a  set and $X^{-1}$ a disjoint set in one-one
correspondence with $X$ via a map $x \rightarrow x^{-1}$ and define
$(x^{-1})^{-1} = x$. We extend this to a map on $(\X)^*$ by defining
$(x_{1}x_{2} \cdots x_{n})^{-1} = x_{n}^{-1} \cdots
x_{2}^{-1}x_{1}^{-1}$, giving $(\X)^*$ the structure of the free
monoid with involution on $X$. Throughout this paper by an {\em
$X$-graph} (or just an {\em edge-labeled graph} if the labeling set
$X$  is understood) we mean a strongly connected digraph $\graph$
with edges labeled over the set $\X$ such that the labeling is
consistent with an involution: that is, there is an edge labeled $x
\in \X$ from vertex $v_1$ to vertex $v_2$ if and only if there is an
inverse edge labeled $x^{-1}$ from $v_2$ to $v_1$. The initial
vertex of an edge $e$ will be denoted by $\init e$ and the terminal vertex by $\ter e$. {The $X$-graph $\graph$
with one vertex and one positively labeled edge labeled by $x$ for each $x \in X$ is referred to as the {\em bouquet of $|X|$ circles} and is denoted by $B_X$.} If $X = \emptyset$, then we view $\graph$ as the
graph with one vertex and no edges.

The label on an edge $e$ is denoted by $\ell(e) \in \X$. There is an
evident notion of {\em path} in an $X$-graph.  The initial (resp.
terminal) vertex of a path $p$ will be denoted by ${\alpha}(p)$
(resp. ${\omega}(p)$). The label on the path $p = e_{1}e_{2} \ldots
e_k$ is the word $\ell(p) = \ell(e_{1})\ell(e_{2})\ldots \ell(e_{k})
\in (\X)^*$.

$X$-graphs occur frequently in the literature.
The Cayley graph ${\graph}(G,X)$ of a group $G$ relative to a set
$X$ of generators is an $X$-graph: its vertices are the elements of
$G$ and it has an edge labeled by $x$ from $g$ to $gx$ for each $x
\in \X$.

 If we designate an initial vertex (state) $\alpha$ and a
terminal vertex (state) $\beta$ of $\graph$, then the birooted
$X$-graph ${\mathcal A} = ({\alpha}, {\graph}, {\beta})$ may be
viewed as an automaton. See for example the book of Hopcroft and
Ullman \cite {HU} for basic information about automata theory. The
language accepted by this automaton is the subset $L({\mathcal A})$
of $(\X)^*$ consisting of the words in $(\X)^*$ that label paths in
$\graph$ starting at $\alpha$ and ending at $\beta$. This automaton
is called an {\em inverse automaton} if it is deterministic (and
hence co-deterministic), i.e. if for each vertex $v$ of $\graph$
there is at most one edge with a given label starting at $v$ or
ending at $v$. This also implies that any path is uniquely
determined by its initial vertex and its label. A deterministic $X$-graph can also be defined as an edge-labeled graph $\Gamma$ that admits a label-preserving graph morphism to $B_X$ which,  for any vertex $v$ of $\Gamma$, is injective on the set of edges $e$ with $\alpha(e) = v$.
Such maps are called {\it graph immersions} in
\cite{stall}.


For each subset $N$ of an inverse monoid $M$, we denote by
$N^{\omega}$ the set of all elements $m \in M$ such that $m \geq n$
for some $n \in N$.  The subset $N$ of $M$ is called {\em closed} if
$N = N^{\omega}$.

Closed inverse submonoids of an inverse monoid $M$ arise naturally
in the representation theory of $M$ by partial injections on a set
 \cite{schein}. An inverse monoid $M$ acts (on the right) by
 injective partial functions on a set $Q$ if there is a homomorphism
 from $M$ to $\SIM(Q)$. Denote by $q.m$ the image of $q$ under the
 action of $m$ if $q$ is in the domain of the action by $m$.
 The following basic fact is well known (see \cite{schein}).

 \begin{Prop}
\label{stabclosed}
If an inverse monoid $M$ acts on $Q$ by injective partial functions, then for every $q \in
 Q, \, \Stab(q) = \{m \in M : q.m = q\}$ is a closed inverse submonoid
 of $M$.
\end{Prop}

 Conversely, given a  closed inverse submonoid $H$ of $M$, we
 can construct a transitive representation of $M$ as follows. A
 subset of $M$ of the form $(Hm)^{\omega}$ where $mm^{-1} \in H$ is
 called a {\em right ${\omega}$-coset} of $H$. Let $X_H$ denote the
 set of right $\omega$-cosets of $H$. If $m \in M$, define an action
 on $X_H$ by $Y . m = (Ym)^{\omega}$ if $(Ym)^{\omega} \in X_H$ and
 undefined otherwise. This defines a transitive action of $M$ on
 $X_H$ with $\Stab(H)=H$. Conversely, if $M$ acts transitively on $Q$, then this action
 is equivalent in the obvious sense to the action of $M$ on the
 right $\omega$-cosets of $\Stab(q)$ in $M$ for any $q \in Q$. See
 \cite{schein} or \cite{Pet} for details.

 The {\em $\omega$-coset graph} ${\graph}_{(H,X)}$ (or just ${\graph}_H$ if $X$ is understood)
 of a closed inverse
 submonoid $H$ of an $X$-generated inverse monoid $M$  is
 constructed as follows. The
 set of vertices of ${\graph}_H$ is $X_H$ and there is an edge
 labeled by $x \in \X$ from $(Ha)^{\omega}$ to $(Hb)^{\omega}$ if
 $(Hb)^{\omega} = (Hax)^{\omega}$.
 Then ${\graph}_H$ is a
 deterministic
 $X$-graph. The
 birooted $X$-graph $(H,{\graph}_{H},H)$ is called the {\em
 $\omega$-coset automaton} of $H$. The language accepted by this
 automaton is $H$ (or more precisely the set of words $w \in (\X)^*$
 whose natural image in $M$ is in $H$).
 Clearly, if $G$ is a group generated by $X$, then ${\graph}_H$ coincides with the coset graph of
 the
 subgroup $H$ of $G$.

We call two closed inverse submonoids $H_1, H_2$ of an inverse
monoid $M$ \textit{conjugate} if there exists $m \in M$ such that
$mH_1m^{-1} \subseteq H_2$ and $m^{-1}H_2m \subseteq H_1$. It is
clear that conjugacy is an equivalence relation on the set of closed
inverse submonoids of $M$: however, conjugate closed inverse
submonoids of an inverse monoid  are not necessarily isomorphic. For
example, the closed inverse submonoids $\{1,aa^{-1},a^{2}a^{-2}\}$
and $\{1,aa^{-1},a^{-1}a, aa^{-2}a\}$ of the free inverse monoid on
the set $\{a\}$ are conjugate by $a^{-1}$ but not isomorphic.

Here we note that since inverse monoids form a variety of algebras
(in the sense of universal algebra --- i.e. an equationally defined
class of algebras), free inverse monoids exist. We will denote the
free inverse monoid on a set $X$ by $\FIM(X)$. This is the quotient
of $(\X)^*$, the free monoid with involution, by the congruence that
identifies $ww^{-1}w$ with $w$ and $ww^{-1}uu^{-1}$ with
$uu^{-1}ww^{-1}$ for all words $u,w \in (\X)^*$. See \cite{Pet} or
\cite{law} for much information about $\FIM(X)$. In particular,
\cite{Pet} and \cite{law} provide an exposition of  Munn's solution
\cite{munn} to the word problem for $\FIM(X)$ via birooted
edge-labeled trees called {\em Munn trees}.

In his thesis \cite{Ste2} and paper \cite{Ste1}, Stephen initiated
the combinatorial theory of presentations of inverse monoids by extending Munn's
results about free inverse monoids to arbitrary presentations of
inverse monoids. We refer the reader to \cite{Ste1} or our paper
\cite{MSz} for details of Stephen's construction of Sch\"utzenberger
graphs and Sch\"utzenberger automata and their use in the study of
presentations of inverse monoids.

We recall that an {\em inverse category} is a category $\cat$ with the property that for every morphism $p$
in $\cat$ there is a unique inverse morphism $p^{-1}$ such that
$p=pp^{-1}p$ and $p^{-1}=p^{-1}pp^{-1}$. Note this implies $\alpha(p)=\omega(p^{-1})$ and $(p^{-1})^{-1}=p$. In this paper, the inverse categories we consider will consists of certain equivalence classes of paths in graphs, thus the domain of $p$ is denoted by $\alpha(p)$, the codomain by $ \omega(p)$, and we multiply morphisms from left to right.
The {\em loop monoids} $L(\cat,
v)$ of an inverse category, that is, the set of all morphisms from
$v$ to $v$, where $v$ is an arbitrary vertex, are inverse monoids.

Analogously to inverse semigroups, a natural partial order can be defined on the morphisms of each inverse category $C$ by putting $p \leq q$ if $p=qe$ for some idempotent $e$. Note that in this case the morphisms $p$ and $q$ are necessarily coterminal.
Each inverse category $\cat$ has a greatest groupoid homomorphic image, obtained by identifying two morphisms $p$ and $q$ if and only if there exists a morphism $t$ in $\cat$ such that $t \leq p$ and $t \leq q$ in the natural partial order on $\cat$, and in this case $p$ and $q$ are both coterminal with $t$ and thus with each other. Alternatively, the morphisms identified are precisely the pairs in the congruence generated by $pp^{-1} \approx 1_{\alpha(p)}$, where $p$ is any morphism.

\begin{Prop}
\label{fundgpoid}
Let $\cat$ be an inverse category and $v$ a vertex in $\cat$. Then the maximal group homomorphic image of $L(\cat , v)$ is isomorphic to the vertex group of $v$ in the maximal groupoid image of $\cat$. 

\end{Prop}

\begin{proof}
Two loops $p$ and $q$ at $v$ in $\cat$ are identified in the maximal groupoid image of $\cat$ if and only if they have a common lower bound $ep=eq$ in the natural partial order, and that lower bound must be in $L(\cat, v)$ since idempotents in $\cat$ are loops. Thus $p$ and $q$ are identified in the maximal groupoid image of $\cat$ if and only if they are identified in the maximal group image of the loop monoid $L(\cat , v)$.

\end{proof}

\subsection{$CW$-complexes and $\Delta$-complexes}

As our results build on previous, analogous results for $CW$-complexes \cite{MSz}, we regard $\Delta$-complexes as $CW$-complexes with special attaching maps, following the definition in \cite{hatch}.
Recall the following definition of a finite dimensional
$CW$-complex $\comp$:
\begin{enumerate}
    \item Start with a discrete set $\comp^0$, the $0$-cells of $\comp$.
    \item Inductively, form the $n$-skeleton $\comp^n$ from $\comp^{n-1}$ by attaching $n$-cells
    $C^n_\tau$ via \emph{attaching maps}
$\varphi_\tau \colon S^{n-1} \to \comp^{n-1}$. This means that
$\comp^n$ is the quotient space of $\comp^{n-1}\ \dot\sqcup_\tau\
B^n_\tau$ under the identifications $x \sim \varphi_\tau (x)$ for $x
\in \bd B_\tau^n$. The cell $C^n_\tau$ is a homeomorphic image of
$B^n_\tau - \bd B^n_\tau$ under the quotient map. (Here $B^n$ is the
unit ball in ${\mathbb R}^n$ and $S^{n-1} = \partial B^n$ is its
boundary).
     \item Stop the inductive process after a finite number of steps to obtain
a finite dimensional $CW$-complex $\comp$.
 \end{enumerate}

  The dimension of the
complex is the largest dimension of one of its cells. Note that a $1$-dimensional $CW$-complex is just
an undirected graph, with the usual topology. We denote the
set of $n$-cells of $\comp$ by $\comp^{(n)}$. Each 
cell $C^n_\tau$ is open in the topology of
 the $CW$-complex $\comp$. A subset $A \subseteq \comp$ is open iff
 $A \cap \comp^n$ is open in $\comp^n$ for each $n$. We emphasize that each cell of the complex is an open cell by definition. 

 Each cell $C^n_\tau$ has a {\it characteristic map}
 ${\sigma}_\tau$, which is defined to be the composition $B^n_\tau
 \hookrightarrow (\comp^{n-1}\ \dot\sqcup_\tau\
B^n_\tau) \rightarrow \comp^n \hookrightarrow \comp$. This is a continuous
map whose restriction to the interior of $B^n_\tau$ is a
homeomorphism onto $C^n_\tau$ (equipped with a subspace topology) and whose restriction to the boundary
of $B^n_\tau$ is the corresponding attaching map $\varphi_\tau$. An
alternative way to describe the topology on $\comp$ is to note that
a subset $A \subseteq \comp$ is open iff ${\sigma}_{\tau}^{-1}(A)$
is open in $B^n_\tau$ for each characteristic map ${\sigma}_\tau$.


The standard $n$-simplex is the set

\begin{center}

${\Delta}^n = \{(t_{0},...,t_{n}) \in {\mathbb R}^{n+1} :
{\Sigma}_{i} t_{i} = 1$ and $t_{i} \geq 0$ for all $i \}$.

\end{center}

We denote the $n+1$ {\it vertices} of ${\Delta}^n$ by  $v_{i} =
(0,\ldots, 0, 1, 0,\ldots , 0)$ ($1$ in $i$th position).  We order
vertices by $v_{i} < v_{j}$ if $i < j$. The {\it faces} of the
simplex are the  subsimplices with vertices any non-empty subset of
the $v_{i}$'s. There are $n+1$ faces of dimension $n-1$, namely the
faces ${\Delta}^{n-1}_i = [v_{0},...,v_{i-1},v_{i+1},...v_{n}]$ for
$i = 0,1,...,n$ spanned by omitting one vertex.

A $\Delta$-complex is a  quotient space of a collection of disjoint
simplices obtained by identifying certain of their faces via the
canonical linear homeomorphisms that preserve the ordering of
vertices.

Equivalently, a $\Delta$-complex is a $CW$-complex $X$ in which each
$n$-cell $e_{\alpha}^n$ has a characteristic map
${\sigma}_{\alpha} : {\Delta}^{n} \rightarrow X$ such that the
restriction of ${\sigma}_{\alpha}$ to each $(n-1)$-dimensional face
of ${\Delta}^n$ is the characteristic map for an
$(n-1)$-cell of $X$. 


The order on the vertices of the simplex makes it naturally possible
to regard each $k$-cell $C^{k}_{\tau}$ of a $\Delta$-complex
$\mathcal C$ as a {\it rooted} cell, with distinguished root the
image under the characteristic map ${\sigma}^{k}_{\tau}$ of the
minimal $0$-cell in the order on $0$-cells in ${\Delta}^k$. We will
denote the root of the cell $C$ by ${\alpha}(C)$. Thus we may regard
the $1$-skeleton as a digraph with each $1$-cell (edge) $e$ directed
from its initial vertex (the root of the cell) to its terminal
vertex $\omega(e)$ (the image of the maximal $0$-cell of
${\Delta}^1$ under the characteristic map). It is also convenient to introduce a {\em ghost edge} $e^{-1}$ for each $1$-cell (edge) $e$  of the $1$-skeleton of $\mathcal C$. Here $\alpha(e^{-1}) = \omega(e), \, \omega(e^{-1}) = \alpha(e)$ and the set $\{e^{-1} : e \in \mathcal{C}^{(1)}\}$ is assumed to be disjoint from and in one-one correspondence with $\mathcal{C}^{(1)}$. We extend the notation by defining $(e^{-1})^{-1} = e$ for each $e \in \mathcal{C}^{(1)}$. For technical reasons,
if $C$ is a cell of at least $2$ dimensions, we define
$\omega(C)=\alpha(C)$.


In the sequel, we will further assume all complexes to be connected.

\subsection{Immersions between $\Delta$-complexes}

We call a map $f \colon \mathcal D \to \mathcal C$ between $\Delta$-complexes a \emph{$\Delta$-map}
if for each $k$-cell $D^k_\tau$ of
$\mathcal D$, $f(D^k_{\tau})$ is a $k$-cell of $\mathcal C$ \colored{and $f$ commutes with the characteristic maps; i.e. $f \circ
{\sigma}^{k}_{\tau} = {\gamma}^{k}_{\tau}$ where ${\sigma}^{k}_{\tau}$ is the characteristic map of $D^{k}_{\tau}$ and ${\gamma}^k_{\tau}$ is the characteristic map of $f(D^{k}_{\tau})$}.
The notion of a $\Delta$-map comes from regarding $\Delta$-complexes as geometric realizations of $\Delta$-sets (also called semi-simplicial sets). Such maps are also known to be continuous \cite{Milnor}. Also note that their restrictions to (open) cells are homeomorphisms onto their images and they map the root of a cell to the root of its image.

For  a vertex $v$ of a $\Delta$-complex $\mathcal C$ we define
$\star_{\mathcal C}(v) = \{e: e$ is an edge or a ghost edge or a cell of dimension $k \geq 2$ in $\mathcal C$ with $\alpha(e) = v\}$.
Then a $\Delta$-map $f : \mathcal D \rightarrow \mathcal C$ induces a map $f_v : \star_{\mathcal D}(v) \rightarrow \star_{\mathcal C}(f(v))$ in the obvious way.

\begin{Def}
\label{immersion}
An \emph{immersion} from a $\Delta$-complex $\mathcal D$ to a $\Delta$-complex
$\mathcal C$ is a $\Delta$-map $f : {\mathcal D} \rightarrow
\mathcal C$ for which each induced map $f_v : \star_{\mathcal D}(v) \rightarrow \star_{\mathcal C} (f(v))$ is injective.
\end{Def}

We remark that for graphs ($1$-dimensional complexes) this concept coincides with the definition of immersion between graphs in the sense of Stallings \cite{stall}. 
In the Appendix, we show that an immersion is precisely a locally injective $\Delta$-map. Hence topological immersions are immersions in the sense of Definition \ref{immersion}, and in the locally compact case, the converse also holds.

\section{Generalized paths}

A \emph{path} in a $\Delta$-complex $\comp$ is a path on its $1$-skeleton in the graph theoretic sense, where we allow edges to be read in either direction. Recall that the $1$-skeleton can be regarded as a digraph since the edges of a simplex are oriented from vertices with smaller indices to vertices with larger indices. Also each edge $e$ comes equipped with a ghost edge $e^{-1}$.
A path then is a
sequence of directed edges and ghost edges
$e_{1}e_{2}...e_n$ with $\omega (e_{i}) = \alpha (e_{i+1})$.
We also consider empty paths $1_v$ around every $0$-cell $v$.

A \emph{generalized path} is a sequence
$e_{1}e_{2}...e_n$ where each $e_i$ is either a $k$-cell with $k \geq
1$ or a ghost edge, and $\omega (e_{i-1})
= \alpha (e_{i})$ for $i = 2,...,n-1$. (Recall that $\omega (e) = \alpha
(e)$ if $e$ is a $k$-cell with $k \geq 2$.) If $p=e_{1}\ldots e_n$, then we define $\alpha(p)=\alpha(e_1)$ and $\omega(p)=\omega(e_n)$. We define $p^{-1}$ to be $e_n^{\epsilon_n} \ldots e_1^{\epsilon_1}$, where
$$\epsilon_i=
\begin{cases}
-1,\hbox{ if }e_i\hbox{ is a $1$-cell or a ghost edge}\\
1\hbox{ otherwise.}
\end{cases}$$
Note that $p^{-1}$ is a generalized path from $\omega(p)$ to $\alpha(p)$, and $(p^{-1})^{-1}=p$. The length of a generalized path 
$p=e_{1}e_{2}...e_n$ is $n$, the length of $p=1_v$ is $0$. We will denote the length of $p$ by $|p|$.

In order to define the equivalence relation on generalized paths that will give rise to loop monoids, we need the notion of boundary paths. If ${\Delta}^{k} = [v_{0},v_{1},...,v_{k}]$, then ${\Delta}^k$ has
$k+1$ faces of dimension  $k-1$, namely the $(k-1)$-simplices
${\Delta}^{k-1}_{i}$ for
$i = 0,...,k$. All of these faces  except ${\Delta}^{k-1}_{0}$
contain the vertex $v_0$. The smallest vertex of
${\Delta}^{k-1}_{0}$ under the order on vertices is $v_1$.

If $\comp$ is a $\Delta$-complex  and $C^{k}$
is a $k$-cell of $\comp$, there is a corresponding characteristic
map ${\sigma}^{k}\colon {\Delta}^{k} \rightarrow {\mathcal C}$. The
restriction of ${\sigma}^{k}$ to ${\Delta}^{k-1}_i$ is a
characteristic map ${\sigma}^{k-1}_{i}$ of some $(k-1)$-dimensional
cell $C^{k-1}_{i}$ of $\mathcal C$, by definition of a
$\Delta$-complex. The root of $C^{k}$ is ${\alpha}(C^{k}) =
{\sigma}^{k}(v_{0})$, and the root of $C^{k-1}_i$ is also
${\sigma}^{k}(v_{0})$ if $i \neq 0$, but the root of $C^{k-1}_0$ is
${\sigma}^{k}(v_{1})$. Thus the $1$-cell
${\sigma}^{k}([v_{0},v_{1}])$ is a directed edge in the $1$-skeleton
of $\mathcal C$ from the root of $ C^k$ to the root of $C^{k-1}_0$.

For a $2$-cell $C^{2}$ of a $\Delta$-complex $\mathcal C$, the \emph{boundary path}\footnote{This is called boundary walk in \cite{MSz}.} $bp(C^2)$ of $C^2$ is the image in $\mathcal C$ of the path
$[v_{0},v_{1}][v_1,v_{2}][v_{0},v_2]^{-1}$ in the $1$-skeleton of ${\Delta}^{2} =
[v_{0},v_{1},v_{2}]$ under the corresponding characteristic map from
${\Delta}^2$ to $\mathcal C$. For a $k$-cell $C^{k}$ of  $\mathcal
C$, denote the image under ${\sigma}^{k}$ of the $1$-cell
$[v_{0},v_{1}]$ of ${\Delta}^k$ by $e(C^{k})$. For $k \geq 3$,
the \emph{boundary path} $bp(C^k)$ of $C^k$ is the generalized path
$C^{k-1}_{k}C^{k-1}_{k-1}\ldots C^{k-1}_{1}e(C^{k})C^{k-1}_{0}(e(C^{k}))^{-1}$.

\subsection{The loop monoid $L(\comp, v)$}
We define a category $\operatorname{FC}(\comp)$ with objects $\comp^0$ and morphisms generalized paths of $\comp$ --- this is the analogue of the free category on a graph. The domain of $p$ is its initial vertex and its codomain its terminal vertex, thus, as mentioned, $\alpha(p)$ and $\omega(p)$ are unambiguous.
Let $\sim$ denote the congruence on $\operatorname{FC}(\comp)$ generated by the relations $p \sim pp^{-1}p$, and $pp^{-1}qq^{-1} \sim qq^{-1}pp^{-1}$ for all
generalized paths with $\init p = \init q$, and the additional relations $CC \sim C$ and $C \sim C\, bp(C)$ for all $C \in \comp^{(k)}, k \geq 2$. Let $\pcat \comp$ denote the category $\operatorname{FC}(\comp)/\sim$. Note that $\pcat \comp$ is an inverse category. We define the \emph{loop monoid} of $\comp$ at $v \in \comp^{(0)}$ to be the inverse monoid $L(\pcat \comp, v)$. We denote it by $L(\comp, v)$. This consists of the $\sim$-classes of generalized paths around $v$.

A $\Delta$-map $f \colon \compb \to \comp$ induces a map $f \colon \operatorname{FC}(\compb) \to \operatorname{FC}(\comp)$ defined on the objects as before, and on the morphisms by $f(e_1 \ldots e_n)=f(e_1)\ldots f(e_n)$ for any non-empty generalized path $e_1 \ldots e_n$, and $f(1_v)=1_{f(v)}$ for $v \in \compb^{(0)}$. Note that the induced map is a functor, i.e. $\alpha(f(p))=f(\alpha(p))$, $\omega(f(p))=f(\omega(p))$, and $f(p)f(q)=f(pq)$ for any connecting generalized paths $p$ and $q$ in $\compb$. Notice furthermore that $f(p)^{-1}=f(p^{-1})$ and $f(bp(C))=bp(f(C))$ holds for any generalized path $p$ and any cell $C$ of dimension greater than $1$. This implies that $f$ respects the relations generating $\sim$, and hence whenever $p \sim q$ we have $f(p)\sim f(q)$. Thus $f$ induces a functor $f \colon \pcat\compb \to \pcat\comp$, a fact we will use frequently.

The following lemma formalizes the statement that when paths lift along immersions, they lift uniquely.

\begin{Lem}
\label{lem:immersionpath}
If $f \colon \compb \to \comp$ is an immersion and $p,q$ are generalized paths in $\compb$ with $f(p)=f(q)$ and $\alpha(p)=\alpha(q)$, then $p=q$.
\end{Lem}

\begin{proof}
Note that $f(p)=f(q)$ implies $|p|=|q|=k$. We prove the lemma by induction on $k$. If $k=0$, then by $\alpha(p)=\alpha(q)$ we have $p=q$. Assume that $k \geq 1$, and put $p=p'e_p$, $q=q'e_q$ where $e_p, e_q$ are cells or ghost edges. Note we have $f(p')=f(q')$ and $\alpha(p')=\alpha(q')$ by the assumption that $f(p)=f(q)$ and $\alpha(p)=\alpha(q)$, and thus by induction $p'=q'$. Therefore $\alpha(e_p)=\omega(p')=\omega(q')=\alpha(e_q)$, and denoting this common vertex by $v$ we have $e_p, e_q \in \star_{\compb}(v)$. It follows from the assumption we have $f(e_p)=f(e_q)$ and so we obtain $e_p=e_q$ from star-injectivity. This proves $p=q$.
\end{proof}

\begin{Prop}
\label{fundgroup} For any vertex $v$ in a (connected) $\Delta$-complex
 $\comp$, the greatest group homomorphic image of $L(\comp, v)$ is
the fundamental group of $\comp$ at $v$.
\end{Prop}

\begin{proof}
Recall that the fundamental groupoid of $\comp$ is obtained by factoring the category of paths (on the 1-skeleton) with the congruence generated by relations of the form $pp^{-1} \approx 1_{\alpha(p)}$ for any path $p$, and $bp(C)\approx 1$ for any $2$-cell $C$.
Observe that this is the greatest groupoid homomorphic image of $\pcat \comp$, that is, it is obtained by factoring  $\pcat \comp$ with the congruence generated by
relations of the form $qq^{-1}\approx 1_{\alpha(q)} $ for
any generalized path $q$ (as this implies $C \approx 1_{\alpha(C)}$ when $C \in \comp^{(k)}, k \geq 2$).
It follows by Proposition \ref{fundgpoid} that the fundamental group $\pi_1(\comp , v)$ is the greatest group homomorphic image of the corresponding loop monoid $L(\comp,
v)$.
\end{proof}

\section{Labeled $\Delta$-complexes}

In this section, we introduce a way to assign labels to all cells of
a $\Delta$-complex in a way that cells with the same root have different labels,
and boundary paths of cells with the same label also have the same label.
This will allow us to identify generalized paths with their initial vertex and label, and to think of elements $L(\comp,v)$ as equivalence classes of labels. Labeled $\Delta$-complexes are higher dimensional analogues of $X$-graphs, and, like $X$-graphs, are defined via an immersion into a $\Delta$-complex with one $0$-cell.

\begin{Lem}
\label{immtoB}  Every  $\Delta$-complex $\mathcal C$  admits an
immersion into a $\Delta$-complex with one $0$-cell.
\end{Lem}

\begin{proof}
If we identify all $0$-cells of $\mathcal C$, then the quotient cell
complex is also a $\Delta$-complex $\mathcal B$ with one $0$-cell.
 The corresponding map $f :
{\mathcal C} \rightarrow \mathcal B$ is an immersion (in fact a topological immersion) since it is
injective on all $k$-cells with $k > 0$.

\end{proof}

Let $\mathcal B$ be a $\Delta$-complex of dimension $n $ with one
$0$-cell. Let  $\{e^{1}_{\rho} : {\rho} \in X\}$ be its set of
$1$-cells, $\{e^{k}_{\rho} : {\rho} \in P_{k}\}$  its set of
$k$-cells for $2 \leq k \leq n$, and let
${\beta}^{k}_{\rho} : {\Delta}^{k} \rightarrow {\mathcal B}$ be the
characteristic map of $e^{k}_{\rho}$ for $k \geq 1$. Here we assume
that the sets $X, P_k$ are all mutually disjoint.
We denote this
$\Delta$-complex $\mathcal B$ by $B(X,P_{2},...,P_{n},
\{{\beta}^{k}_{\rho}\})$, or more briefly by $B(X,P)$ where $P =
P_{2} \cup ... \cup P_n$.
We view $X$ as a
set of labels for the $1$-cells of $B(X,P)$ and $P_k$ as a set of
labels of the $k$-cells of $B(X,P)$ for $2 \leq k \leq n$. That is,
the label on the $k$-cell $e^{k}_{\rho}$ is ${\ell}(e^{k}_{\rho}) =
\rho$.

If $\mathcal B$ is one dimensional, we obtain a  bouquet of circles,
usually denoted by $B_X$. One should think of $B(X,P)$ as a
generalization. The $1$-skeleton of $B(X,P)$ is $B_X$, which we
regard as an {\em $X$-graph} as usual; i.e. each edge labeled by $x
\in X$ is equipped with an inverse (ghost) edge labeled by $x^{-1}$. The
labeling on the $1$-cells of $B(X,P)$ extends to a labeling on {\it
paths} in the $1$-skeleton of $B(X,P)$ in the obvious way. The label
on a path $p$ in $B(X,P)$ will be denoted by ${\ell}(p)$: thus
${\ell}(p) \in (\X)^*$. More generally we may extend the labeling on
cells to a labeling on {generalized paths} in the obvious way:
if $e_{1}e_{2}...e_{t}$ is a generalized path, then
${\ell}(e_{1}e_{2}...e_{t}) =
{\ell}(e_{1}){\ell}(e_{2})...{\ell}(e_{t}) \in (X \cup X^{-1} \cup
P)^*$. The label of an empty path is the empty word.

\begin{Def}
Let  $f : {\mathcal C} \rightarrow
B(X,P)$ be an immersion between $\Delta$-complexes and define the label of a  cell $C^{k}_{\tau} \in \mathcal C^{(k)}$ with $k \geq 1$  by
${\ell}(C^{k}_{\tau}) = {\ell}(f(C^{k}_{\tau}))$. We say that  
$\mathcal C$ is {\it labeled} over the
complex $B(X,P)$ via the immersion $f$. 
\end{Def}

Note that
${\ell}(C^{1}_{\tau}) \in X$  and ${\ell}(C^{k}_{\tau}) \in P_{k}$
if $2 \leq k \leq n$, and cells of $\mathcal C$ have the same label
if and only if they have the same image under $f$. If the underlying complex
$B(X,P)$ is understood, we just say that $\mathcal C$ is a {\it
labeled complex}. In particular every complex $B(X,P)$ is labeled via the identity map on itself.

By Lemma \ref{immtoB}, every $\Delta$-complex admits some labeling.
The immersion $f$ constructed in the proof of that lemma assigns
different labels to all cells of $\mathcal C$. Of course we would
usually choose smaller sets $X, P_k$ as sets of labels for the cells
of $\mathcal C$ if possible.

\begin{Rem}
\label{welldef}
Notice that  if $p$ and $q$ are generalized path of a labeled $\Delta$-complex $\mathcal C$ with $\alpha(p)=\alpha(q)$ and $\ell(p)=\ell(q)$, by Lemma \ref{lem:immersionpath} applied to the immersion inducing the labeling, we have $p=q$. In particular if
$C^{k}_{\gamma}$ and $C^{s}_{\tau}$ are
distinct cells of a labeled $\Delta$-complex $\mathcal C$ with the
same root $v$ then ${\ell}(C^{k}_{\gamma}) \neq {\ell}(C^{s}_{\tau})$.
Furthermore, if $e_1$ and $e_2$ are distinct $1$-cells of a $\Delta$-complex with the same terminal vertex, then $\ell(e_1) \neq \ell(e_2)$.
\end{Rem}

\begin{Def}
An immersion $g : {\mathcal D}
\rightarrow \mathcal C$ of $\Delta$-complexes is said to commute with the labeling (or to respect the labeling) if $\comp$ and $\compb$ are labeled over the same complex $B(X,P)$
by immersions $f_\comp \colon \comp \to B(X,P)$ and $f_\compb \colon \compb \to B(X,P)$, and $g$ commutes with these labeling maps, that is, $f_\comp \circ g= f_\compb$.
\end{Def}

Note that if $g : {\mathcal D}
\rightarrow \mathcal C$ is an immersion of $\Delta$-complexes, then a labeling $f_\comp$ 
 of $\mathcal C$ over $B(X,P)$ induces a labeling $f_\compb = f_\comp \circ g $  of
$\mathcal D$ over $B(X,P)$ such that $g$
 respects the labeling.

For any cell $C \in \comp^{(k)}, k\geq 2$, define the boundary label $bl(C)$ of $C$ as $\ell(bp(C))$. We then have the following lemma:

\begin{Lem}\label{boundarylabel}
Let $g : {\mathcal D} \rightarrow \mathcal C$ be an immersion of a
$\Delta$-complexes that commutes with the labeling and let $D^{k}$ be a
$k$-cell of $\mathcal D$ with $k \geq 1$.  Then

\begin{enumerate}
\item ${\ell}(D^{k})={\ell}(g(D^{k}))$,
\item $g({\alpha}(D^{k})) = {\alpha}(g(D^{k}))$, and if $k \geq 2$, then
\item $bl(D^{k}) = bl(g(D^{k}))$.
\end{enumerate}
\noindent In particular, if $f$ is an immersion of $\mathcal D$ into
$B(X,P)$ that defines a labeling of $\mathcal D$, then $bl(D^{k}) =
bl(f(D^{k}))$ for every $k$-cell $D^k$ of $\mathcal D$ where $k \geq 2$.
\end{Lem}

\begin{proof}
The statement is  immediate from the definition and the fact that immersions commute with characteristic maps.
\end{proof}

For the remainder of this paper, we will assume that all
$\Delta$-complexes are labeled and that all immersions between
$\Delta$-complexes respect the labeling, as described above.

\section{Properties of the loop monoid}
\label{mxp}

Note that it follows from Lemma \ref{boundarylabel} that if $\comp$ is a labeled $\Delta$-complex, then any two $k$-cells $(k \geq 2)$ that have the same label also have the same
boundary label. This allows us to define relations on $\FIM(X \cup P)$ analogous to those defining the inverse category $\pcat \comp$. 
For $\rho \in P_k$ we define $bl({\rho})$ to be the common boundary label of all the cells with label $\rho$.

\begin{Def}
Consider a $\Delta$-complex $B(X,P)$. We define $M(X,P)$ as the inverse monoid with
{\it generators} $X \cup P$ and {\it relations}

\begin{itemize}

\item{${\rho}^{2} = \rho$ for each ${\rho} \in P$  and}

\item{${\rho} = {\rho} \, bl({\rho})$
 for each ${\rho} \in P$}.
\end{itemize}
\end{Def}

We remark that the conditions ${\rho} = {\rho}^2$ and ${\rho} =
{\rho} \, bl({\rho})$ are equivalent to ${\rho} = {\rho}^2$ and
 ${\rho} \leq  bl({\rho})$. It is sometimes more convenient to use
 this characterization of the defining relations for
 $M(X,P)$.

\begin{Prop}
\label{propmxp}
The inverse monoid $M(X,P)$ is the loop monoid of $B(X,P)$ based at its unique $0$-cell $v$.

\end{Prop}

\begin{proof}
Observe that the inverse monoid $M(X,P)$ can be defined as the quotient of the free monoid $(X \cup X^{-1}\cup P \cup P^{-1})^\ast$ by the congruence generated by pairs of the form $(w, ww^{-1}w)$ and $(ww^{-1}uu^{-1},uu^{-1}ww^{-1})$ where $u,w$ are arbitrary words, and the additional relations ${\rho}^{2} = \rho$ and ${\rho} = {\rho} \, bl({\rho})$ for each $\rho \in P$. Note that these relations imply $\rho=\rho^{-1}$, therefore $P^{-1}$ can be omitted from the generating set.

In $B(X,P)$, generalized paths are exactly words in $X \cup X ^{-1}\cup P$, and the relation $\sim$ on these words is generated by exactly the above pairs. Hence $L(B(X,P),v)$ is indeed equal to $M(X,P)$.
\end{proof}

 The following two lemmas are crucial to prove that immersions into a $\Delta$-complex are encoded by its loop monoids.

\begin{Lem}
\label{path} 
Let $C$ be any $k$-cell with $k \geq 2$ in a
$\Delta$-complex labeled by an immersion into $B(X,P)$. If $p$ is
a path on ${\partial}C$ obtained as the image under the attaching map of a
 closed path around $v_0$ in the $1$-skeleton of $\Delta^k$,
 then ${\ell}(C) \leq {\ell}(p)$ in $M(X,P)$.
\end{Lem}

\begin{proof}
Let $q$ be a path on
$\Delta^k$ around $v_0$, and let $p$ be its image under the attaching map.
We work by induction on $|p|$. If $|p|=0$, then $\ell(p)=1 \geq \ell(C)$ holds as $\ell(C)$ is an idempotent. Otherwise put $q=f_1 \cdots f_t$ and $p =
e_{1}e_{2} \cdots e_{t}$. Note we cannot have $|p|=1$. If $|p|=2,3$, then $q$ lies in a $2$-face $F$ of $\Delta^k$ containing $v_0$. Let $D$ be the cell corresponding to $F$ under the characteristic map of $C$. We claim that $\ell(D) \geq \ell(C)$. This can be seen by induction on $k$: if $k=2$, then $D=C$ so it is clear; if $k >2$, let $\Delta_i$ be a $(k-1)$-face of $C$ containing $F$ and hence $v_0$. For the corresponding cell $C_i$, we have $\ell(C_i) \geq bl(C) \geq \ell(C)$ by the presentation as $bl(C)$ is the product of $\ell(C_i)$ and some idempotents, and $\ell(D) \geq \ell(C_i)$ by induction.

Suppose $p=|3|$, then we must have $q=f_1f_2f_3$ where $bp(D)=(f_1f_2f_3)^{\pm1}$. 
Then by the presentation, $\ell(e_1e_2e_3), \ell(e_1e_2e_3)^{-1}\geq \ell(D) \geq \ell(C)$, which gives the statement $\ell(p) \geq \ell(C)$.
If $|p|=2$, then retaining the notation of the previous case, we must have $q$ of the form $f_1f_1^{-1}$ with $bp(D)=(f_1f_2f_3)^{\pm1}$. Then $\ell(p)=\ell(e_1e_1^{-1}) \geq \ell(e_1e_2e_3)\ell(e_1e_2e_3)^{-1} \geq \ell(C)$ as well.

Suppose $|p| > 3$. 
Take a path of minimal length from ${\omega}(f_{2})$ to ${\alpha}(f_{1}) = v_0$ in
$\Delta^k$ (that is, of length $0$ or $1$),
and denote its image under the characteristic map of $C$ by $s_1$. Consider the path 
$(e_1e_2s_1)(s_1^{-1} e_3 \ldots e_t)$, and note that both factors satisfy the assumption of the lemma, but have length $< |p|$. Thus by induction, we obtain $\ell(p)\geq \ell(e_1e_2s_1)\ell(s_1^{-1} e_3 \ldots e_t) \geq \ell(C)$.
\end{proof}

We extend this to obtain a similar result about
generalized paths.

\begin{Lem}
\label{genpath} Let $C$ be any $k$-cell with $k \geq 2$ in a
$\Delta$-complex labeled over $B(X,P)$. If $p$ is the image  under
the corresponding attaching map of any
 closed generalized path around $v_0$ on the boundary of $\Delta^k$,
 then ${\ell}(C) \leq {\ell}(p)$ in $M(X,P)$.

\end{Lem}
\begin{proof}
Consider a closed generalized path
$q=q_{1}\Delta_{1}q_{2}\Delta_{2}...q_{t}\Delta_{t}q_{t+1}$ around $v_0$ on $\bd
\Delta^k$, where $q_i$ is a path, and $\Delta_i$ is a face of
$\Delta^k$ of at least $2$-dimensions, and let
$p=p_{1}D_{1}p_{2}D_{2}...p_{t}D_{t}p_{t+1}$ be the image of $q$  under the
attaching map. Let $s_i$ be the image of the shortest path from $v_0$ to $\alpha(\Delta_i)$ on $\bd \Delta^k$ (this is either an edge or the empty path).
We claim that $\ell (C) \leq \ell (s_i)\ell(D_i) \ell(s_i^{-1})$. This can be proved by induction on $n=k-\dim(D_i)$. For $n = 0$ we have $C=D_i$ and $s_i$ the empty path, so the claim is immediate. If $n\geq 1$, then let $\Delta'$ be a $\dim(D_i)+1$-dimensional face of $\Delta^k$ containing $v_0$, and let $D=\sigma_C(\Delta')$. Then $bl(D)$ is the product of $\ell(s_i)\ell(D_i)\ell(s_i^{-1})$ and some idempotents, so $\ell(D) \leq \bl(D) \leq \ell (s_i)\ell(D_i) \ell(s_i^{-1})$ from the presentation, and $\ell(C) \leq \ell(D)$ by induction.

By Lemma \ref{path} we have ${\ell}(C) \leq
{\ell}(p_{1}s_{1}^{-1}), {\ell}(C) \leq
{\ell}(s_{1}p_{2}s_{2}^{-1}),\ldots, {\ell}(C) \leq
{\ell}(s_{t}p_{t+1})$. Hence,  by multiplying all of these
inequalities, we obtain
\begin{align*}
{\ell}(C) & \leq
{\ell}(p_{1}s_{1}^{-1}){\ell}(s_{1}){\ell}(D_{1}){\ell}(s_{1}^{-1}){\ell}(s_{1}p_{2}s_{2}^{-1})
{\ell}(s_{2}){\ell}(D_{2}){\ell}(s_{3}^{-1})\cdots\\
&{\ell}(s_{t-1}p_{t}s_{t}^{-1}){\ell}(s_{t})
{\ell}(D_{t}){\ell}(p_{t+1}) \leq
{\ell}(p_{1}){\ell}(D_{1}){\ell}(p_{2}){\ell}(D_{2})\cdots{\ell}(D_{t}){\ell}(p_{t+1})
= {\ell}(p)
\end{align*}
 as required.
\end{proof}

A small, but useful observation:

\begin{Prop}
\label{Mproperty}
The inverse submonoid of $M(X,P)$ generated by $X$ is isomorphic to $\FIM(X)$.
\end{Prop}

\begin{proof}
Recall that $M(X,P)$ is obtained as
a factor of the free monoid on $X \cup X^{-1} \cup P$ by the defining
relations of inverse monoids and those introduced in the
presentation; these are all equations which, if they contain a
letter in $P$, then they contain it on both sides. Therefore, to a
word in $(X \cup X^{-1})^\ast$, one can only apply those relations
which define the free inverse monoid on $X$.
\end{proof}

\subsection{$L(\comp, v)$ as a stabilizer}

Now let $\mathcal C$ be a $\Delta$-complex labeled over $B(X,P)$.
Then we may define a natural {\it action} of the inverse monoid
$M(X,P)$  by partial one-to-one maps on the set ${\mathcal C}^{(0)}$
of $0$-cells of $\mathcal C$ as follows. For $x \in \X$ and $v \in
{\mathcal C}^{(0)}$ define $v.x = w$ if there is an edge (or a ghost edge) labeled by
$x$ from $v$ to $w$ in the $1$-skeleton of $\mathcal C$,
 and $v.x$ is undefined if there is no such edge. For
${\rho} \in P_{k}$ with $k \geq 2$ and $v \in {\mathcal C}^{(0)}$,
define $v.{\rho} = v$ if $v = {\alpha}(C^k)$ for some $k$-cell
$C^{k}$
 with ${\alpha}(C^{k}) = v, \, \ell (C^k) = \rho$, and $v. {\rho}$ is undefined otherwise.

We begin by a simple, but useful lemma about the above action.

\begin{Lem}
\label{l:Bensnewlemma}
 Let $\comp$ be a $\Delta$-complex labeled over $B(X,P)$.
\begin{enumerate}
\item For any generalized path $p$ in $\comp$, $\alpha(p).\ell(p)$ is defined and is equal to $\omega(p)$.
\item For any $w \in (X \cup X^{-1} \cup P)^\ast$, if $v.w$ is defined for some $v \in \comp^{(0)}$, then there exists a generalized path in $\comp$ with $\alpha(p)=v$, $\ell(p)=w$, and $\omega(p)=v.w$.
\end{enumerate}
\end{Lem}

\begin{proof}
\begin{enumerate}
\item We prove the statement by induction on the length of $p$. If $p$ is an empty path, then $\ell(p)=1$ and $\omega(p)=\alpha(p)$, so the statement holds. Otherwise put $p=p'e$ with $e$ a cell or a ghost edge. Then by induction $v.\ell(p')$ is defined and is equal to $\omega(p')$, and notice that $\alpha(e)=\omega(p')$. Then by definition we have $\omega(p').\ell(e)=\omega(e)=\omega(p)$, and so $v.\ell(p)=v.\ell(p').\ell(e)=\omega(p)$ as needed. 
\item The proof is by induction on the length of $w$. If $w$ is the empty word, then $v.w=v$ and $p=1_v$ suffices. Otherwise put $w=w'y$ with $y \in X \cup X^{-1} \cup P$. As $v.w=v.w'.y$, we have that $v.w'$ is defined and thus there exists by induction a generalized path $p'$ in $\comp$ with $\alpha(p')=v$, $\ell(p')=w'$, $\omega(p')=v.w'$. Then $v.w=\omega(p').y$, thus $\omega(p').y$ is defined. Since $y \in  X \cup X^{-1} \cup P$, by the definition of the action, there is a cell or ghost edge $e$ with $\alpha(e)=\omega(p')$, $\ell(e)=y$,  $\omega(e)=v.w$. Put $p=p'e$, this is indeed a generalized path with $\alpha(p)=v$, $\omega(p)=v.w$ and $\ell(p)=\ell(p')\ell(e)=w'y=w$.
\end{enumerate}
\end{proof}

\begin{Lem}
\label{action}
The above action of the generators $X \cup X^{-1} \cup P$ of $M(X,P)$ on ${\mathcal
C}^{(0)}$ extends to a well-defined action of $M(X,P)$ on ${\mathcal
C}^{(0)}$.
\end{Lem}

\begin{proof}
By the Remark \ref{welldef}, the labeling we introduced ensures that letters in $X \cup X^{-1}\cup P$
act by partial
one-to-one maps on ${\mathcal C}^{(0)}$, and the action of $x \in X$ is the inverse of the action of $x^{-1}$. Thus the action of the
generators extends to a well-defined action of the free inverse
monoid  $\FIM(X \cup P)$ on ${\mathcal C}^{(0)}$. We need to show
that the action respects the defining relations of $M(X,P)$. The
 actions of letters in $P$ are, by definition, identity maps
on their domains, hence idempotent maps. 
If $\rho \in P_2$ fixes $v$ and $\ell (C^2) = \rho$ with $\alpha (C^2) = v$ then necessarily $bp(C^2)$ determines a closed path at $v$ and so $bl(\rho)$ fixes $v$. So
the
action by $\rho$ is the same as the action by ${\rho} \, bl({\rho})$
for ${\rho} \in P_{2}$. Finally, for ${\rho} \in P_{k}, k \geq 3$,
the action by $\rho$ is also a restriction of the action by
$bl({\rho})$. This is because each of the elements
${\ell}(C^{k-1}_{i})$ (for $i = 1,..., k$) and
${\ell}(e(C^{k})C^{k-1}_{0}(e(C^{k}))^{-1})$ that arise in the
definition of the boundary label $bl(C^k)$ of a cell $C^k$ stabilize
${\alpha}(C^k)$. So again the action by $\rho$ is the same as the
action by ${\rho} \, bl({\rho})$. Hence the action respects all
defining relations in $M(X,P)$, as required.
\end{proof}

The stabilizer $\Stab(v)$ of a $0$-cell $v \in {\mathcal C}^{(0)}$ under the action by
$M(X,P)$ is a closed inverse submonoid of $M(X,P)$, and stabilizers
of different vertices in ${\mathcal C}^{(0)}$ are conjugate closed
inverse submonoids of $M(X,P)$. Hence the immersion $f: {\mathcal C}
\rightarrow B(X,P)$ that defines the labeling of $\mathcal C$ gives
rise to a conjugacy class of closed inverse submonoids of $M(X,P)$.

\begin{Prop}
For any labeled $\Delta$-complex $\comp$ and any $v \in \comp^{(0)}$, we have $\Stab(v) \cong L(\comp, v)$.
\end{Prop}

\begin{proof}
We show that the map between $L(\comp, v)$ and  $\Stab(v)$ taking the $\sim$-class of a generalized path $p$ to the value of $\ell(p)$ in $M(X,P)$ is an isomorphism.

Denote the $\sim$-class of a generalized path $p$ in $\comp$ by $[p]$ and the equivalence class of a word $w \in (X \cup X^{-1} \cup P)^*$ under the congruence $\approx$ on $(X \cup X^{-1} \cup P)^*$ defining the inverse monoid $M(X,P)$ by $[[w]]$. 
We claim that for any pair of coterminal generalized paths $p,q$, we have $p \sim q$ if and only if $\ell(p) \approx \ell(q)$.
Note that for $w_1,w_2 \in (X \cup X^{-1} \cup P)^*, \, w_1 \approx w_2$ if and only if it is possible to pass from $w_1$ to $w_2$ by a sequence of elementary transitions involving the relations used to define the free inverse monoid on $X \cup P$ together with the additional relations used to define $M(X,P)$ as a quotient of this free inverse monoid. Similarly, $p \sim q$ if it is possible to pass from $p$ to $q$ by a sequence of elementary transitions involving paths, using the relations defining the category $\pcat \comp$. It is clear that any sequence of such elementary transitions of paths defines a sequence of elementary transitions of their labels, so $p \sim q$ implies $\ell(p) \approx \ell(q)$.

For the converse, let $p$ be any generalized path, and $w$ a word in $X^\ast$ such that $\ell(p)$ and $w$ are related by an elementary transition. We claim that then $w$ labels a generalized path $q$ coterminal with $p$, and hence there is a corresponding path transition from $p$ to $q$. If the transition from $\ell(p)$ to $w$ uses an identity defining inverse monoids, then this is easy to see. The same is true if $w$ uses a relation of the form $\ell(C)^2=\ell(C)$ as $\ell(C)$ always labels a closed path. For the last type of relation $\ell(C)\bl(C)=\ell(C)$, notice that whenever we have a (closed) path labeled by $\ell(C)$ at some vertex $u$, we also have a closed path labeled by $\bl(C)$ at $u$. Thus we can freely remove or add $\bl(C)$ after $\ell(C)$ and obtain a generalized path. 

This implies that for any coterminal paths $p,q$, if we have any sequence of elementary transitions from $\ell(p)$ to $\ell(q)$, there must be a sequence of corresponding elementary path transitions from $p$ to a path which is coterminal with $p$ and has label $\ell(q)$, but this must be $q$. Thus $\ell(p) \approx \ell(q)$ implies $p \sim q$.

We claim that the map $\phi \colon [p] \mapsto [[\ell(p)]]$ is an isomorphism from from $L(\comp, v)$ to $\Stab (v)$. First note that by the previous assertion, for any pair of closed generalized paths $p,q$ around $v$, we have $[p]=[q]$ if and only if $[[\ell(p)]]=[[\ell(q)]]$. This shows that $\phi$ is injective, and a well-defined homomorphism into $M(X,P)$.
Furthermore $[p] \in L(\comp, v)$ if and only if $p$ is a closed generalized path in $\comp$ at $v$, which by Lemma \ref{l:Bensnewlemma} implies $v.\ell(p)=v$, that is $[[\ell(p)]] \in \Stab(v)$. Conversely if $[[w]] \in \Stab(v)$ then by Lemma \ref{l:Bensnewlemma} there is a generalized path $p$ with $\ell(p)=w$ and $v.\ell(p)=v$, thus $\phi([p])=[[w]]$. This shows that the image of $\phi$ is $\Stab(v)$ indeed. This completes the proof.
\end{proof}

In the sequel, we identify $\sim$-classes of $(v,v)$-paths with their common label in $M(X,P)$, and regard $L(\comp, v)$ equal to $\Stab(v)$.

\section{Classification of immersions}
\label{mainsec}

The previous section has shown how an immersion into $B(X,P)$ gives
rise to a conjugacy class of closed inverse submonoids of its
(unique) loop monoid $M(X,P)$. This section shows that this is true
for any immersion, moreover, the converse also holds: immersions are
in one-to-one correspondence with conjugacy classes of closed
inverse submonoids of loop monoids.

\begin{Prop}
\label{imm1}
Let $\comp$ and $\compb$ be $\Delta$-complexes labeled over a common complex $B(X,P)$, and
suppose $g \colon \compb \to \comp$ is an immersion, and let $v \in \compb^{(0)}$. Then
$L(\compb, v)$ is a closed inverse submonoid of $L(\comp, f(v))$.
\end{Prop}

\begin{proof}
There are immersions $f_\comp \colon \comp \to B(X,P)$ and $f_\compb \colon \compb \to B(X,P)$
defining the labeling, and $f_\comp \circ g=f_\compb$. Then $L(\compb, v)=\Stab(\compb,v)$ and
$L(\comp, f(v))=\Stab(\comp, f(v))$ are both closed inverse submonoids of the corresponding $M(X,P)$. Therefore it is enough to show that $\Stab(\compb,v) \subseteq \Stab(\comp, f(v))$. Indeed, suppose $w$ labels a closed path $p$ around $v$ in $\compb$. Then $g(p)$ is a closed path around $f(v)$ with the same label $w$, which proves the statement.
\end{proof}

We proceed to develop the theory of the  converse part of the
correspondence.
Fix a $\Delta$-complex $\comp$ labeled over some
$B(X,P)$ by an immersion $g \colon \comp  \to B(X,P)$.
We define a
graph $\graph_\comp$ associated with $\comp$ as
follows:

\begin{center}

$V(\graph_\comp)= \comp^{(0)}$ and

$E(\graph_\comp)= \comp^{(1)} \cup \{f_{C^{k}}: C^{k}$ is a $k$-cell
in $\comp^{(k)}, k \geq 2\},$

\end{center}

\noindent where $f_{C^{k}}$ denotes a loop based at $\alpha (C^k)$ and
labeled by $\ell (C^{k})$.  Thus the edges
 in $\comp^{(1)}$ are labeled over $\X$ and the edges of the form $f_{C^{k}}$ (for $C^k$ a
 $k$-cell)
 are
 labeled over $P_{k}$.
Since an edge labeled by $\rho \in P_k$ is always a loop labeled by
an idempotent in $M(X,P)$, we may identify $P_k$ with
 $P_{k}^{-1}$ and regard $\graph_\comp$ as  an $(X \cup P
 )$-graph in the sense defined in the introduction.
 If we abuse notation slightly by identifying the loop $f_{C^{k}}$ in $\graph_\comp$
 with the $k$-cell $C^k$,
 then  paths in
 the graph $\graph_\comp$ are just identified with generalized paths in the
 $\Delta$-complex $\mathcal C$.

 Let $u \in
\comp^{(0)}$. Given a closed inverse submonoid $H$ of $L(\comp, u)$,
we construct a complex $\comp_H$ that immerses into $\comp$ with $H
=L(\comp_H, \tilde u)$ for some $\tilde u \in \comp_H^{(0)}$ \colored{in the} preimage of $u$.

The complex $\comp_H$ is defined with the help of the $\omega$-coset graph $\graph_H$ of $H$.
The inverse monoid
$\M{X,P}$ acts on the vertices of
$\graph_H$. The idempotents, of course, all label loops.

We build a $\Delta$-complex $\mathcal C_H$  such that $\Gamma_{\mathcal
C_H}=\Gamma_H$. This complex has the following sets of cells:
$$\comp_H^{(0)}=V(\graph_H),$$
$$\comp_H^{(1)}=\{e \in E(\graph_H): \ell (e) \in X\},$$
$$\hbox{if }k\geq 2,\ \comp_H^{(k)}=\{\tilde C_e: e \in E(\graph_H) \hbox{ and } \ell (e) \in P_k\}.$$
The attaching maps of $1$-cells are the ones inherited from the graph $\Gamma_H$.
Note that the $1$-skeleton then immerses into the $1$-skeleton of $\comp$ via some map
$f_{1} \colon \comp_H^{1} \to \comp^{1}$ by the results of \cite{MM1},
as Proposition \ref{Mproperty} ensures that our definitions and assumptions
reduce to those of \cite{MM1} when applied to graphs. Moreover, the vertex $H$ is mapped to $u$.

We build the rest of the attaching maps of $\comp_H$ inductively.
Let $k \geq 2$ and denote the set of edges of $\Gamma_H$ with labels in $P_k$ by $E_k$. Suppose $\comp_H^{k-1}$ is a $\Delta$-complex with $1$-skeleton $\comp_H^1$
such that there exists an immersion $f_{k-1}\colon \comp_H^{k-1} \to
\comp$ with $H \mapsto u$, and $\Gamma_{ \comp_H^{k-1}}$ is the subgraph of $\Gamma_H$
induced by the edges labeled by $X \cup P_2 \cup \ldots \cup P_{k-1}$. Note we then have $V(\Gamma_{ \comp_H^{k-1}})=V(\Gamma_H)$, because $\Gamma_H$ is connected and $P_2 \cup \ldots \cup P_{k-1}$ all label loops.

Let $e \in E_k$, then $e$ is a loop. Put $\rho=\ell(e)$ and $v=\alpha(e)=\omega(e)$.
Let $p$ be a path in $\Gamma_H$ (and so in $\comp_H^{k-1}$) from $H$ to $v$. Then $\ell(p)\rho\ell(p)^{-1} \in H \subseteq L(\comp, u)$, and so $\rho$ labels a $k$-cell $C_e$ in $\comp$ based at $\omega(f_{k-1}(p))=f_{k-1}(v)$. Denote its attaching map by $\varphi_e$.

\begin{Lem}
\label{liftattach}
There exists a unique attaching map $\tilde\varphi_e :\partial\Delta^k \to\comp_H^{k-1}$ of a
$k$-cell $\tilde{C}_e$ based at $v$ such that $f_{k-1} \circ \tilde\varphi_e=\varphi_e$.
\end{Lem}

\begin{proof}
We can write the boundary of a $k$-simplex as a (non-disjoint) union of $(k-1)$-dimensional simplices:
$$\partial\Delta^k=\bigcup_{i=1}^{k+1}\Delta_i^{k-1}.$$

By the definition of a $\Delta$-complex, $\varphi_e=\bigcup_{i=1}^{k+1}\sigma_i^{k-1}$ for the
characteristic maps $\sigma_i^{k-1}\colon \Delta_i^{k-1} \to \comp$ of some cells $C_i^{k-1}$, and
likewise, if such a $\tilde\varphi_e$ exists, it has to be the form of
$\tilde\varphi_e=\bigcup_{i=1}^{k+1}\tilde\sigma_i^{k-1}$ for the characteristic maps
$\tilde\sigma_i^{k-1}: \Delta_i^{k-1} \to \comp_H^{k-1}$ of some cells $\tilde C_i^{k-1}$.
Since $f_{k-1} \circ \tilde\varphi_e=\varphi_e$, we have that $\tilde C_i^{k-1}$ is a preimage of
$C_i^{k-1}$ under $f_{k-1}$, in particular, $\ell(C_i^{k-1})=\ell(\tilde C_i^{k-1})$.

For every $i$ between $1$ and $k+1$, let $p_i$ be a path on the one-dimensional faces (edges) of $\Delta^k$ from $v_0$ to the root of $\Delta_i^{k-1}$. Put $q_i=\varphi_e(p_i)$, these are paths on the one-cells of $\comp$. 
Since $p_i \Delta_i^{k-1}(\Delta_i^{k-1})^{-1} p_i^{-1}$ is a closed generalized path on $\partial\Delta^k$ around $v_0$, for the path $\varphi_e(p_i \Delta_i^{k-1}(\Delta_i^{k-1})^{-1} p_i^{-1})=q_i C_i^{k-1} (C_i^{k-1})^{-1}q_i^{-1}$ in $\comp$ around $f_{k-1}(v)$, by Lemma \ref{genpath} we have
$$\rho=\ell(C^k)\leq\ell(q_i)\ell(C_i^{k-1})\ell(C_i^{k-1})^{-1}\ell(q_i)^{-1},$$
therefore $\ell(q_i)\ell(C_i^{k-1})\ell(C_i^{k-1})^{-1}\ell(q_i)^{-1}$ labels a closed path in $\Gamma_H$  around $v$, and therefore in $\comp_H^{k-1}$ also. 
The cell $\tilde {C_i}^{k-1}$ thus can only be the unique $(k-1)$-cell with the label $\ell(C_i^{k-1})$ occurring in the previous path by Remark \ref{welldef}, and $\tilde{\sigma_i}^{k-1}$ is the characteristic map corresponding to $\tilde C_i^{k-1}$. 
Note that for this cell we have $f_{k-1}(\alpha(\tilde C_i^{k-1}))=\alpha(C_i^{k-1})$ and so
$f_{k-1}(\tilde C_i^{k-1})=C_i^{k-1}$ indeed, thus $f_{k-1} \circ \tilde\sigma_i^{k-1}=\sigma_i^{k-1}$.

It remains to be shown that the map $\tilde\varphi_e=\bigcup_{i=1}^{k+1}\tilde\sigma_i^{k-1}$, given as a union of maps on non-disjoint domains, is well-defined;
that is, for any intersection $\Delta_{i,j}=\Delta_i \cap \Delta_j$, we have
$\tilde{\sigma_i}^{k-1}|_{\Delta_{i,j}}=\tilde{\sigma_j}^{k-1}|_{\Delta_{i,j}}$.
Since $\Delta_{i,j}$ is a face of both $\Delta_i$ and $\Delta_j$, both maps
$\tilde{\sigma_i}^{k-1}|_{\Delta_{i,j}}$ and $\tilde{\sigma_j}^{k-1}|_{\Delta_{i,j}}$
are characteristic maps $\tilde\sigma_i^{k-2}$ and $\tilde\sigma_j^{k-2}$ for some
$k-2$-cells $\tilde{C_i}^{k-2}$ and $\tilde{C_j}^{k-2}$. Of course, since
$\varphi_e=\bigcup_{i=1}^{k+1}\sigma_i^{k-1}$, we have
${\sigma_i}^{k-1}|_{\Delta_{i,j}}={\sigma_j}^{k-1}|_{\Delta_{i,j}}=\varphi_k|_{\Delta_{i,j}}$,
hence

\begin{equation}
\label{eqn}
f_{k-1} \circ \tilde\sigma_i^{k-2}={\sigma_i}^{k-1}|_{\Delta_{i,j}}={\sigma_j}^{k-1}|_{\Delta_{i,j}}=f_{k-1} \circ \tilde\sigma_j^{k-2},
\end{equation}

Now take a path $s_i$ on the edges of $\Delta_i$ from the root of
$\Delta_i$ to the root of $\Delta_{i,j}$, and similarly a path $s_j$
on $\Delta_j$. Let $t_i=\sigma_i^{k-1}(s_i)$,
$t_j=\sigma_j^{k-1}(s_j)$, and $\tilde t_i=\tilde
\sigma_i^{k-1}(s_i)$, $\tilde t_j=\tilde \sigma_j^{k-1}(s_j)$. Since
$p_i s_i s_j^{-1} p_j^{-1}$ is a closed path on $\Delta^k$ around
$v_0$,  for the path $\varphi(p_i s_i s_j^{-1} p_j^{-1})=q_it_i
t_j^{-1} q_j^{-1}$ in $\comp$, we have $\ell(q_it_i t_j^{-1}
q_j^{-1})\geq \rho$ by Lemma \ref{genpath}, and there is a closed path at $v$ labeled by
$\ell(q_it_i t_j^{-1} q_j^{-1})$ in $\comp_H^{k-1}$. But, the unique
path from $v$ labeled by $\ell(q_i)$ ends in
$\alpha(\tilde{C_i}^{k-1})$ --- that is how $\tilde{C_i}^{k-1}$ was
defined ---, and the unique path from $\alpha(\tilde{C_i}^{k-1})$
labeled by $\ell(t_i)$ is $\tilde{t_i}$. The same can be said about
$\alpha(\tilde{C_j}^{k-1})$ and $\tilde{t_j}$, which implies that
$\ter{\tilde{t_i}}=\ter{\tilde{t_j}}$, hence $\alpha(\tilde
{C_i}^{k-2})=\alpha(\tilde{C_j}^{k-2})$. Since $f_{k-1}$ is an
immersion which maps $\tilde{C_i}^{k-2}$ and $\tilde{C_j}^{k-2}$ to
the same cell by (\ref{eqn}), this immediately implies
$\tilde{C_i}^{k-2}=\tilde{C_j}^{k-2}:=\tilde{C}^{k-2}$, and hence
$\tilde\sigma_i^{k-2}=\tilde\sigma_j^{k-2}$. That proves that
$\tilde\varphi_e$ is well-defined, and by nature of the
construction, unique. Note that $\tilde\varphi_e$ is continuous,
since it is a union of continuous maps defined on closed sets,
completing the proof.

\end{proof}

We now define $\comp^k_H$ as the $\Delta$-complex with $(k-1)$-skeleton $\comp^{k-1}_H$ and $k$-cells $\comp_H^{(k)}=\{\tilde C_e \colon e \in E_k\}$,  where the attaching map of a cell $\tilde C_e$ is $\tilde\varphi_e$, as defined in the previous lemma. Notice that we by construction have that 
$\Gamma_{\comp_H^k}$ is the subgraph of $\Gamma_H$ induced by the edges labeled by $X \cup P_2 \cup \ldots \cup P_{k}$.

\begin{Lem}
\label{liftattach2}
There is a unique immersion $f_k\colon {\comp^k_H} \to \comp$ for which $f_k|_{\comp_H^{k-1}}=f_{k-1}$.
\end{Lem}

\begin{proof}
Let $\tilde{\sigma}_e$ denote the characteristic map of the cell $\tilde{C}_e$ of $\comp^k_H$ for any $e \in E_k$, and let ${\sigma}_e$ denote the characteristic map of respective cell $C_e$ of $\comp$.
Let $f_k \colon \comp^k_H \to \comp$ be the map defined by
\begin{equation}
  \label{ftul}
  \begin{aligned}
    f_k|_{\comp_H^{k-1}}&=f_{k-1}\\
    f_k|_{\tilde{C_e}^k}&=\sigma_e|_{\int \Delta^k} \circ (\tilde \sigma_e|_{\int \Delta^k})^{-1}
  \end{aligned}
\end{equation}
for $e \in E_k$.
We show that $f_k$ is an immersion: it suffices to show that $f_k$ commutes with the characteristic
maps, and induces injections between star sets. 

Since $f_{k-1}$ is an immersion, $f_k$ clearly commutes with the
characteristic map of any cell contained in ${\comp_H^{k-1}}$. To see
that $f_k$ commutes with the characteristic map of a $k$-cell $\tilde{C}_e$ where $e \in E_k$, let
$x \in \Delta^k$ be an arbitrary point, and consider $f_k \circ
\tilde{\sigma}_e(x)$. If $x \in \int \Delta^k$, then $f_k \circ
\tilde{\sigma}_e(x)= \sigma_e(x)$ by (\ref{ftul}). If $x \in
\bd\Delta^k$, then $x$ lies in a simplex $\Delta^j$ on
$\bd\Delta^k$, and $\tilde{\sigma}_e|_{\Delta^j}$ is a
characteristic map \colored{that commutes with $f_k$}, therefore we again obtain $f_k
\circ \tilde{\sigma}_e(x)= \sigma_e(x)$, as desired.

%

Also, since by construction $\Gamma_{\comp_H^k}$ is a subgraph of $\Gamma_H$, thus whenever $C_1, C_2$ are cells of $C_H^k$ with $\ell(C_1)=\ell(C_2)$ and $\alpha(C_1)=\alpha(C_2)$, we have $C_1=C_2$. In particular if $C_1, C_2 \in \star({\mathcal{C}}^k_H,v)$ for some $v$, then $f(C_1)=f(C_2)$ implies $C_1=C_2$, thus 
 $f$ is an immersion.

The uniqueness of $f$ follows from the fact that any map satisfying
the conditions of the lemma must satisfy (\ref{ftul}).
\end{proof}

We are ready to state and prove the main theorems of the paper. If $\comp$ is a $\Delta$-complex with $u \in \comp$, we call the pair $(\comp,u)$ a pointed $\Delta$-complex. A pointed $\Delta$-map between the pointed $\Delta$-complexes $(\comp,u) \to (\compb,v)$ is a $\Delta$-map $f \colon \compb \to \comp$ with $f(u)=v$. In particular we can talk about pointed immersions and pointed isomorphisms.

Our first main result says that the pointed immersions into a $\Delta$-complex are up to pointed isomorphisms classified by the closed inverse submonoids of its loop monoid.

\begin{Thm}
\label{imm2} Let $\comp$ be a $\Delta$-complex  labeled over some
$B(X,P)$, let $u \in \comp^{(0)}$, and let $H$ be  any closed
inverse submonoid of $L(\comp, u)$. Then  there exists a 
pointed $\Delta$-complex $(\comp_H,v)$ and a pointed immersion $f \colon (\comp_H,v) \to (\comp,u)$ with $H =L(\comp_H, v)$. Furthermore, if $(\comp_H',v')$ is another pointed $\Delta$-complex with an immersion $f'\colon (\comp_H',v') \to (\comp,u)$, $H= L(\comp_H', v')$, then there is a pointed isomorphism $g \colon (\comp_H',v') \to (\comp_H, v)$ such that $f'=f \circ g$.
\end{Thm}

\begin{proof}
The existence part of the theorem is clear from the previous
construction. For uniqueness up to pointed isomorphism, note that $H =L(\comp_H, v)=L(\comp_H', v')$ dictates
that there is a pointed isomorphism between both $(\Gamma_{\comp_H}, v)$,  $(\Gamma_{\comp'_H}, v')$ and the $\omega$-coset graph $(\Gamma_H, H)$, which in particular implies there is a pointed isomorphism $g_1$ between their subgraphs induced by the $X$-labeled edges, and hence between $(\comp_H^1, v)$ and $(\comp'^1_H, v')$. Note that if $p$ is an path in $\comp'^1_H$ with $\alpha(p)=v'$, then $\alpha(g_1(p))=v$ and hence $f'(p)$ and $f \circ g_1(p)$ are both paths in $\comp$ starting at $f(v)=f'(v')=u$ with label $\ell(p)$, thus they are equal. As $\comp'^1_H$ is connected, any of its $1$-cells lies on such a path $p$, which shows $f \circ g_1=f'$ .

It follows by induction from the uniqueness of the construction in Lemma
\ref{liftattach} that $g_1$ extends to an isomorphism $g \colon (\comp_H',v') \to (\comp_H, v)$. Indeed, notice that if $g_{k-1} \colon (\comp_H^{k-1}, v) \to (\comp'^{k-1}_H, v')$ is an isomorphism with $f \circ g_{k-1}=f'$ and $e \in E_k$, 
denoting the constructed cells of 
$\comp_H^{k}$, $\comp_H'^{k}$ by $\tilde C_e$, $\tilde C_e'$ with attaching maps $\tilde\varphi_e$, $\tilde\varphi'_e$ respectively,
then $g_{k-1} \circ \tilde\varphi'_e$ and $\tilde\varphi_e$ both satisfy the conditions of Lemma \ref{liftattach} and are thus equal. Hence putting $g_k(\tilde C_e')=\tilde C_e$ for $e \in E_k$, $g|_{\comp_H'^{k-1}}=g_{k-1}$ defines an isomorphism $g_k \colon (\comp'^k_H,v') \to (\comp^k_H, v)$ with $f'=f \circ g_k$.
\end{proof}

An immediate consequence of Proposition \ref{imm1} and
Theorem \ref{imm2} is our first main theorem, characterizing pointed, connected immersions over finite-dimensional, connected $\Delta$-complexes.

\begin{Thm}
\label{main} Let $(\comp,u)$ and $(\compb,v)$ be connected
$\Delta$-complexes labeled over a common $\Delta$-complex $B(X,P)$,
and suppose $f \colon (\compb,v) \to (\comp,u)$ is a connected, pointed immersion that
commutes with the labeling maps, then $f$ induces an embedding of $L(\compb, v)$
into $L(\comp, u)$. 

Conversely, let $\comp$ be a $\Delta$-complex
labeled over a $\Delta$-complex $B(X,P)$, and let $H$ be a closed
inverse submonoid of the corresponding inverse monoid $\M{X,P}$ such that $H
\subseteq L(\comp, u)$ for some $u \in \comp^0$. Then there exists a
pointed $\Delta$-complex $(\compb,v)$ labeled over the same $B(X,P)$, with $L(\compb, v)=H$, and there is a pointed
immersion $f \colon (\compb,v) \to (\comp,u)$.

Furthermore, if $(\compb',v')$ is another pointed $\Delta$-complex with an immersion $f'\colon (\compb',v') \to (\comp,u)$, $H= L(\compb', v')$, then there is a pointed isomorphism $g \colon (\compb',v') \to (\compb, v)$ such that $f'=f \circ g$.
\end{Thm}

The following theorem states that immersions into a complex $\comp$ up to isomorphism are characterized by the conjugacy classes of closed inverse submonoids of any of its loop monoids.

\begin{Thm}
\label{conjugate}
Take a pointed $\Delta$-complex $(\comp, u)$ labeled over $B(X,P)$.
If and $(\compb, v)$, $(\compb',v')$ are pointed $\Delta$-complexes with pointed immersions $f,f'$ into $(\comp, u)$ respectively, then there is an isomorphism $g \colon \compb' \to \compb$ with $f'=f \circ g$ if and only if
$L(\compb, v)$ is conjugate to $L(\compb', v')$ in $L(\comp,u)$.
\end{Thm}

\begin{proof}
First, suppose such an isomorphism $g$ exists. 
Let $m$ label a path from $v$ to $g(v')$ in $\compb$, then $m$
labels a path from $f(v)=u$ to $f(g(v'))=f'(v')=u$ in $\comp$, hence $m \in
L(\comp, u)$. If $w$ labels a closed generalized path in $\compb'$ at $v'$, then it also labels a closed generalized path at $g(v')$ in $\mathcal D$, so $mwm^{-1}$ labels a closed generalized path in $\compb$ at $v$, and $m L(\compb', v') m^{-1}
\subseteq L(\compb, v)$. Similarly $m^{-1} L(\compb, v) m \subseteq
L(\compb', v')$, hence they are conjugate in $L(\comp, u)$ indeed.

For the converse, suppose  $L(\compb, v)$ is conjugate to $L(\compb', v')$ in $L(\comp, u)$. Then there
exists some $m \in L(\comp, u)$ such that 
$$m^{-1}L(\compb, v)m\subseteq L(\compb', v') \hbox{ and } mL(\compb', v')m^{-1} \subseteq L(\compb, v),$$ in particular, $mm^{-1} \in L(\compb,
v)$. Therefore $m$ labels a generalized path from $v$ to some $0$-cell $z$ in $\compb$, and note that $f(z)=u$ by Remark \ref{welldef} since $f(v)=u$ and $m \in L(\comp, u)$. 
We will show that $L(\compb', v')=L(\compb, z)$.
If
$k \in L(\compb', v')$,  then $mkm^{-1} \in L(\compb, v)$ so $mkm^{-1}$ labels a closed generalized path around $v$ in $\compb$, hence $k$
labels closed generalized path around $z$. Thus we have $L(\compb', v') \subseteq L(\compb, z)$.
 On the other hand, if $n \in L(\compb, z)$, then $mnm^{-1}$ labels a closed generalized path around $v$, so $mnm^{-1}
\in L(\compb, v)$, and $m^{-1}mnm^{-1}m \subseteq L(\compb', v')$. Since $L(\compb', v')$ is
closed and $m^{-1}mnm^{-1}m \leq n$, this yields $n \in L(\compb', v')$,
therefore $L(\compb', v')=L(\compb, z)$. The existence of $g$ then follows from the last statement of Theorem \ref{main} applied to $(\compb', v')$ and $(\compb, z)$.

\end{proof}

\section{Closing remarks}

We remark that the constructions of the inverse monoid $M(X,P)$ and
of the complex associated with a closed inverse submonoid of
$M(X,P)$ are effective. The proof of the following theorem makes use
of Stephen's construction of Sch\"utzenberger graphs \cite{Ste1} and
an extension of this developed in \cite{MSz}.  We note that the
result is somewhat surprising in view of the fact that the maximal
group image of $M(X,P)$ is the fundamental group of $B(X,P)$, which
may have undecidable word problem. However, the fact that $M(X,P)$
is not $E$-unitary enables $M(X,P)$ to have decidable word problem
while its maximal group image may not necessarily have decidable
word problem. The proof follows closely along the lines of the proof
of Theorem 5.7 of \cite{MSz}, so we will omit it.

\begin{Thm}
\label{mainalg} (a) If $X$ and $P$ are finite sets, then the word
problem for $M(X, P)$ is decidable.

(b) If $X$ and $P$ are finite sets and $H$ is a finitely generated
closed inverse submonoid of $M(X, P)$, then the  associated
$\Delta$-complex $\comp_H$ is finite and effectively constructible.

\end{Thm}

Similarly, one may obtain the following characterization of the
covering maps. Again the proof closely follows the proof of Theorem
6.1 of \cite{MSz}.

\begin{Thm}
\label{cover} Let $\comp, \compb$ be $\Delta$-complexes labeled by
an immersion over some complex $B(X,P)$, let $f \colon \comp \to
\compb$ be an immersion that respects the labeling, and let $v \in
\comp^0$ be an arbitrary $0$-cell. Then $f$ is a covering map if and
only if $L(\comp, v)$ is a full closed inverse submonoid of
$L(\compb, f(v))$, that is, it contains all idempotents of
$L(\compb, f(v))$.
\end{Thm}

We conclude by raising the question as to whether an extension of
some of the ideas contained in this paper may be developed to
provide a classification of immersions between more general
topological spaces (for example for arbitrary $CW$-complexes). It
would also be of interest to provide a ``presentation-free"
characterization of the inverse category $\pcat \comp$ that serves
the role of the fundamental groupoid in covering space theory.

\section{Appendix}

In this section we discuss the relationship between our definition of immersions between $\Delta$-complexes and the more classical concept of topological immersions. Recall that a $\Delta$-map $f \colon \compb \to \comp$ is locally injective if every  point in $\compb$ has a neighborhood $U$ with $f|_U \colon U \to f(U)$ injective. It is a topological immersion if we furthermore require that $f|_U$ is a homeomorphism between $U$ and $f(U)$, where both sets are equipped with the subspace topology.

We begin with a rather technical, but useful lemma.
\begin{Lem}
\label{l:locinjtechnical}
Let $f \colon \compb \to \comp$ be a $\Delta$-map, and $u$ be any point in $\compb$. Then $f$ is locally injective at $u$ if and only if for any cells $D_1, D_2$ of $\compb$ with $f(D_1)=f(D_2)$ and  $\sigma_{D_1}^{-1}(u) \cap \sigma_{D_2}^{-1}(u) \neq \emptyset$, we have $D_1=D_2$.
\end{Lem}

\begin{proof} Let $D_1, D_2$ be cells as in the statement. Notice that $f(D_1)=f(D_2)$ implies $D_1$, $D_2$ are of the same dimension $n$. Assume $D_1 \neq D_2$ for contradiction.
Let $N$ be an arbitrary neighborhood of $u$ in $\mathcal D$, and take the set
$U:=\sigma_{D_1}^{-1}(N) \cap \sigma_{D_2}^{-1}(N) \subseteq \Delta^n$. This is an open set in $\Delta^n$, and it is also nonempty, as it contains $\sigma_{D_1}^{-1}(u) \cap \sigma_{D_2}^{-1}(u)$.
Take a point $x \in U \setminus \bd
\Delta^n$ -- there certainly exists such a point, as $U$ is nonempty and
open --, and let $x_i = \sigma_{D_i}(x)$. Then as $x_i
\in D_i$, we have $x_1 \neq x_2$, and $f \circ \sigma_{D_1}=f \circ
\sigma_{D_2}$ implies $f(x_1)=f(x_2)$. As $x_i \in N$, we obtain that $f$ is not locally injective at $u$.

For the converse, assume that $f$ is not locally injective at $u$; we need to show that cells $D_1$ and $D_2$ exist as above. 
Note we must have $\dim(\compb) \geq 1$ for local injectivity to fail. Let $D$ be the unique cell containing $u$, and let $\sigma_D \colon \Delta^l \to \ol D$ be its characteristic map. (If $u$ is a $0$-cell, its characteristic map just maps the $0$-simplex to $u$.) Consider the preimage $\sigma_D^{-1}(u)$ of $u$, which is a point in $\int (\Delta^l)$. Regard simplices as metric spaces with the Euclidean metric they inherit from $\mathbb R^n$, and if $l >0$, let
$$M=d(\sigma_D^{-1}(u), \partial\Delta^l).$$
Notice $M >0$. If $l=0$, put $M=0$. For the $l$-simplex $S$, define
$$S^M=\{x \in S: d(x, \partial S)\geq M\}.$$
This is a closed subset of $\int (S)$, containing $\sigma_D^{-1}(u)$.
For any $k \in \mathbb N$, let
$$m_k=\min\{d(F^M,H): F, H \hbox{ are faces of }\Delta^k \hbox{ with }\dim(F)=l, F \not\subseteq H\},$$
note $m_k \in \mathbb R^+ \cup \{\infty\}$.
Of course if $k \geq 1, l$, then $\Delta^k$ has an $l$-dimensional face $F$, for which $F^M \neq \emptyset$, and a face $H$ with $F \not\subseteq H$, so in this case $m_k \in \mathbb R^+$.
Let
$$m=\min \{m_k : k \leq \dim(\compb)\},$$
where $\dim{\compb}$ denotes the dimension of $\compb$. Notice that $m \in \mathbb R^+$.

\textbf{An important observation.}
Let  $C$ be any $k$-cell of $\compb$ whose closure contains $u$, and denote its characteristic map by $\sigma_C$. For any $x \in \sigma_C^{-1}(u)$, let $S_x$ be the unique face of the simplex $\Delta^ k$ containing $x$ in its interior. Then $S_x$ is the $l$-dimensional simplex $\Delta^l$ and $\sigma_C|_{S_x}=\sigma_D$, so $x \in \sigma^{-1}_{D}(u) \subseteq S_x^M$. Suppose $d(x,y)< m$ for some $y \in \Delta^k$, and let $H$ be any face of $\Delta^k$ with $y \in H$. Then of course $d(S_x^M,H)<m \leq m_k$, so the definition of $m_k$ implies $S_x \subseteq H$.

Let $N_u$ be the neighborhood of $u$ defined as follows. For each cell $C$ of $\mathcal D$, consider its characteristic map $\sigma_C: \Delta^k \to \compb$, and let $U_C$ be the open subset of $\Delta^k$ defined by
$$U_C=\bigcup\{B_{m/2}(x): x \in \sigma_C^{-1}(u)\},$$
where $B_{m/2}(x)$ is the open ball of radius $m/2$ at $x$ in $\Delta^k$.
Notice that if $u \notin \overline{C}$, then $U_C =\emptyset$.

Let $$N_u = \bigcup \{\sigma_C(U_C) : C \hbox{ is a cell of } \compb\}.$$
Clearly $u \in N_u$.
We also claim that for any $k$-cell $C$ in $\compb$, $\sigma^{-1}_C(N_u)=U_C$. The containment $\sigma^{-1}_C(N_u)\supseteq U_C$ is immediate from the definition.

 We need to show that $\sigma^{-1}_C(N_u)\subseteq U_C$. Assume $z \in \sigma^{-1}_C(N_u)$, and let $S_z$ be the unique face of $\Delta^k$ containing $z$ in its interior. Then $\sigma_C|_{S_z}$ is a characteristic map $\sigma_B$ of some cell $B$.
Let $y = \sigma_B(z) \in B \cap N_u$. Then $y \in \sigma_{C'}(U_{C'})$ for some $k'$-cell $C'$ with $B \subseteq \ol C'$.
We can rewrite $y \in \sigma_{C'}(U_{C'})$ as $d(\sigma_{C'}^{-1}(y), \sigma_{C'}^{-1}(u)) < \frac{m}{2}$, that is $d(\tilde{y}, \tilde{u}) < \frac{m}{2}$ for some $\tilde{y}\in \sigma_{C'}^{-1}(y)$, $\tilde{u} \in \sigma_{C'}^{-1}(u)$. Let $S_u$ be the unique face of $\Delta^{k'}$ containing $\tilde{u}$ in its interior, and $S_y$ the unique face containing $\tilde{y}$ in its interior. Then by the important observation, $d(\tilde{y}, \tilde{u}) < \frac{m}{2}$ yields that $S_u \subseteq S_y$, in particular $\tilde{u} \in S_y$. Notice that $\sigma_{C'}(S_y)=\overline{B}$, so $\sigma_{C'}|_{S_y}=\sigma_{B}$, and $\sigma_{C'}|_{S_y}^{-1}(y)=\sigma_{B}^{-1}(y)=\{z\}$, hence $\tilde y=z$.
By $\sigma_C|_{S_z}=\sigma_B$ we also have
$\tilde u \in \sigma_{C'}|_{S_y}^{-1}(u)=\sigma_{B}^{-1}(u)=\sigma_{C}|_{S_z}^{-1}(u)$, so
$d(z,\sigma|_C^{-1}(u))\leq d(z,\tilde u)<\frac m 2$, so $z \in U_C$ indeed.

This shows that $\sigma^{-1}_C(N_u)=U_C$, which also immediately yields that $N_u$ is open. It also shows that if $C$ is a cell such that $u \notin\overline{C}$, then $N_u \cap \overline C=\emptyset$.

Let $w_1$ and $w_2$ be distinct points of $N_u$ with $f(w_1)=f(w_2)$. For each $w_i$, there is exactly one cell $D_i$
such that $w_i \in D_i$. Since $f$ restricted to any cell of
$\mathcal D$ is a homeomorphism, $f(w_1)=f(w_2)$ implies $D_1 \neq
D_2$, and $f(D_1)=f(D_2)$; in particular, $D_1$ and $D_2$ must be
of the same dimension, say $n$. Notice that
since $u \in \ol{D_1} \cap \ol{D_2}$,
$D \subseteq \ol D_1 \cap \ol D_2$.

Denote the characteristic map of $D_i$ by $\sigma_{D_i}$, the common cell $f(D_1)=f(D_2)$ by $C$,
the characteristic map of $C$ by $\sigma_{C}$.
Since $\sigma_{D_i}$ restricted to $\int(\Delta^n)$ is a homeomorphism, $\sigma^{-1}_{D_i}(w_i)$ is a single point $\tilde{w_i}$, which satisfies
$$\sigma_{C}(\tilde w_1)=f \circ \sigma_{D_1}(\tilde w_1)=f \circ \sigma_{D_2}(\tilde w_2)=\sigma_{C}(\tilde w_2).$$
But as $\sigma_{C}$ restricted to $\int(\Delta^n)$ is again a homeomorphism, the above implies $\tilde w_1=\tilde w_2$, which we will denote by just $\tilde w$.

Now consider $\sigma_{D_i}^{-1}(N_u)=U_{D_i}$.
Notice that $\tilde{w} \in U_{D_i}$, so there exists $\tilde u_i$ with $\sigma_{D_i}(\tilde u_i)=u$, and $d(\tilde w, \tilde u_i)<\frac{m}{2}$.
We have $\tilde u_i \in \partial\Delta^n$ as $u \in \bd D_i$. Let $F_i$ be the face of $\Delta^n$ with $\tilde u_i \in \int (F_i)$. Then
$$d(\tilde u_1, F_2) \leq d(\tilde u_1 , \tilde w)+d(\tilde w, \tilde u_2 )<m,$$
so by our important observation, we have $F_1 \subseteq F_2$, that is, $F_1=F_2=:F$.
From
$\sigma_{D_i}(F)=D$ we obtain $\sigma_{D_1}|_F=\sigma_D=\sigma_{D_2}|_F$. Denote the preimage of $u$ under the common map $\sigma_{D_i}|_F$ by $\tilde u$. Then $\tilde u \in \sigma^{-1}_{D_1}(u) \cap \sigma^{-1}_{D_2}(u)$, thus $D_1$ and $D_2$ are such that $f(D_1)=f(D_2)$, $\sigma^{-1}_{D_1}(u) \cap \sigma^{-1}_{D_2}(u)\neq \emptyset$, $D_1 \neq D_2$, as required.
\end{proof}

\begin{Prop}
\label{p:immlocinj}
A $\Delta$-map $f \colon \compb \to \comp$ is an immersion if and only if it is locally injective at each $0$-cell of $\compb$.
\end{Prop}

\begin{proof}
Assume first that $f$ is locally injective at each $0$-cell, and let $v \in \compb^0$. Let $D_1, D_2 \in \star_\compb(v)$ and assume $f(D_1)=f(D_2)$. If $D_1, D_2$ are cells, then by $v=\alpha(D_1)=\alpha(D_2)$ we have $v_0 \in \sigma^{-1}_{D_1}(v) \cap \sigma^{-1}_{D_2}(v)$; if they are ghost edges with $D_i=e_i^{-1}$, then similarly $v_1 \in \sigma^{-1}_{e_1}(v) \cap \sigma^{-1}_{e_2}(v)$. In both cases, applying Lemma \ref{l:locinjtechnical}, we obtain $D_1=D_2$. Hence $f$ is injective on the star set of $v$ for any $v \in \compb^0$, thus it is an immersion.

For the converse, assume $f$ is an immersion, and assume for contradiction that $f$ is not locally injective at some $0$-cell $v$. Then by Lemma \ref{l:locinjtechnical} there exist distinct cells $D_1, D_2$ with $f(D_1)=f(D_2)$ and  $\sigma_{D_1}^{-1}(v) \cap \sigma_{D_2}^{-1}(v) \neq \emptyset$. Let $D_1, D_2$ be such examples of minimal dimension $n$. If $\alpha(D_0)=\alpha(D_1)$, then $D_1, D_2 \in \star_\compb(v)$ and so $f$ is not an immersion, a contradiction. So assume this is not the case.

Take $\tilde v \in  \sigma_{D_1}^{-1}(v) \cap \sigma_{D_2}^{-1}(v)$, then $\tilde v \neq v_0$.
If $n=1$, that is, $D_1$ and $D_2$ are edges, then $\tilde v=v_1$. In this case $D_1^{-1}, D_2^{-1} \in \star_\compb(v)$, and so $f$ is not injective on $\star_\compb(v)$, a contradiction.

Now assume $n \geq 2$. As   $f \circ \sigma_{D_1}=f \circ \sigma_{D_2}$,
note that for each proper face $F$ of $\Delta^k$ with $\tilde v \in F$, the cells $D_1'=\sigma_{D_1}(F)$ and $D_2'=\sigma_{D_2}(F)$ also satisfy 
$$f(D_1')=f \circ \sigma_{D_1}(F)=f \circ \sigma_{D_2}(F)=f(D_2'),$$ and $\tilde v \in \sigma_{D_1'}^{-1}(v) \cap \sigma_{D_2'}^{-1}(v)$, hence $D_1'=D_2'$ by the minimality of $n$. Take $F$ to be the edge $[v_0, \tilde v]$ of $\Delta^n$. Then the common cell $\sigma_{D_i}(F)$ is an edge of $\compb$ from $\alpha(D_i)$ to $v$, so $\alpha(D_1)=\alpha(D_2)$. Denoting the latter vertex by $u$, we obtain that $D_1, D_2 \in \star_{\compb}(u)$, which is, again, a contradiction.
\end{proof}

\begin{Prop}
\label{p:locinj2}
A $\Delta$-map $f \colon \compb \to \comp$ is locally injective at $0$-cells if and only if it is locally injective.
\end{Prop}

\begin{proof}
For the nontrivial implication,  let $u \in \mathcal D$ be any point. Take any cells $D_1, D_2$ with $f(D_1)=f(D_2)$, $\sigma_{D_1}^{-1}(u) \cap \sigma_{D_2}^{-1}(u) \neq \emptyset$; by Lemma \ref{l:locinjtechnical} we only need to show that $D_1=D_2$. Let $\tilde u \in \sigma_{D_1}^{-1}(u) \cap \sigma_{D_2}^{-1}(u) \neq 0$, and let $F$ be the unique face of $\Delta^n$ containing $\tilde u$ in its interior. Then $\sigma_{D_i}|_F$ is the characteristic map of the unique cell containing $\sigma_{D_i}(\tilde u)=u$, so $\sigma_{D_1}|_F=\sigma_{D_2}|_F$. Let $\tilde{v}$ be an arbitrary $0$-cell in $F$, and let $v$ be the common $0$-cell $\sigma_{D_i}(\tilde v)$. By assumption, $f$ is locally injective at $v$, and as  $\tilde v \in \sigma_{D_1}^{-1}(v) \cap \sigma_{D_2}^{-1}(v)$, applying Lemma  \ref{l:locinjtechnical} we obtain that $D_1=D_2$ indeed. 
\end{proof}

\begin{Thm}
A $\Delta$-map $f \colon \compb \to \comp$ is an immersion if and only if it is locally injective. Furthermore, if $\comp$ and $\compb$ are locally compact, then $f$ is an immersion if and only if it is a topological immersion.
\end{Thm}

\begin{proof}
The first statement is an immediate consequence of Propositions \ref{p:immlocinj} and \ref{p:locinj2}. In the second statement,
it is clear topological immersions are locally injective and hence immersions. 
For the converse, let $u$ in $\mathcal D$ be any point, by local injectivity there is a neighborhood
$N_u$ of $u$ such that $f|_{N_u} \colon N_u \to f(N_u)$ is a
continuous bijection. As $\mathcal D$ is locally compact, there exists a compact neighborhood $K_u$ of $u$. Let $N$ be an open neighborhood of $u$ such that $\ol N \subseteq K_u \cap N_u$. Then
$f|_{\ol N} \colon \ol N \to f(\ol N)$ is a continuous bijection with $\ol N$  compact and $f(\ol N)$ Hausdorff, therefore it is a homeomorphism. Consequently, $f|_{N} \colon  N \to f(N)$ is also a homeomorphism. This proves that $f$ is a topological immersion.
\end{proof}

\medskip

\noindent {\bf Acknowledgement} The authors are grateful to Benjamin Steinberg for  helpful comments on an earlier version of this paper and to the referee for pointing several errors and misprints in an earlier version: this led to a  revised version of the definition of an immersion between $\Delta$-complexes that is used in this paper.

\bigskip

\end{document}